\documentclass[a4paper,leqno]{amsart}

\usepackage{amsmath,amssymb,amsfonts,epsfig,euscript} 

\def\dis
{\displaystyle}
\def\eps
{\varepsilon}
\def\qed
{\hfill $\square$}

\def\C{{\mathbb C}}
\def\R{{\mathbb R}}
\def\N{{\mathbb N}}
\def\Z{{\mathbb Z}}

\def\F{{\mathcal F}}
\def\Sch{{\mathcal S}}
\def\Eq#1#2{\mathop{\sim}\limits_{#1\rightarrow#2}}
\def\Tend#1#2{\mathop{\longrightarrow}\limits_{#1\rightarrow#2}}

\def\Id{\operatorname{Id}}
\def\d{{\partial}}

\newtheorem{lem}{Lemma}[section]
\newtheorem{cor}[lem]{Corollary}
\newtheorem{prop}[lem]{Proposition}
\newtheorem{theo}[lem]{Theorem}
\newtheorem{defin}[lem]{Definition}
\newtheorem{hyp}[lem]{Assumption}

\def\rq
{\noindent {\it Remark.} }
\def\rqs
{\noindent {\it Remarks.} }

\numberwithin{equation}{section}

\begin{document}
\title[Semi-classical NLS with harmonic potential]{Semi-classical
Schr\"odinger equations with harmonic potential and nonlinear
perturbation} 
\author[R. Carles]{R\'emi Carles}
\address[Former address]{Universit\'e Bordeaux 1\\ Institut
  Math\'ematique de Bordeaux, UMR CNRS 5251\\
351~cours de la Lib\'eration\\
33~405~Talence~cedex\\ France}
\address[Future address, starting September 2007]{Universit\'e
  Montpellier~2\\ Math\'ematiques, 
  UMR CNRS 5149\\ CC 051\\ 
   Place Eug\`ene Bataillon\\ 34095
  Montpellier cedex 5\\ France}
\email{Remi.Carles@math.cnrs.fr}
\begin{abstract}
Solutions of semi-classical
Schr\"odinger equation with isotropic harmonic potential focus
periodically in 
time. We study the perturbation of this equation by a nonlinear
term. If the scaling 
of this perturbation is critical, each focus crossing is described by
a nonlinear scattering operator, which is therefore iterated as many
times as the solution passes through a focus.  The study of this
nonlinear problem is made 
possible by the introduction of two operators well adapted to
Schr\"odinger equations with harmonic potential, and by suitable
Strichartz inequalities.
\\

\noindent\textsc{R\'esum\'e.} Les solutions de l'\'equation de Schr\"odinger
semi-classique avec 
potentiel harmonique isotrope focalisent p\'eriodiquement en
temps. Nous \'etudions la perturbation de cette \'equation par un
terme non-lin\'eaire. Pour une \'echelle critique de cette
perturbation, chaque travers\'ee de foyer est d\'ecrite par un
op\'erateur de diffusion non-lin\'eaire, qui est par cons\'equent
it\'er\'e autant de fois que la solution traverse une caustique. 
Cette
\'etude est permise par l'usage de deux op\'erateurs qui s'av\`erent
bien adapt\'es \`a l'\'equation de Schr\"odinger avec potentiel
harmonique, et par des estimations de
Strichartz ad\'equates. 
\end{abstract}
\maketitle

\section{Introduction}
\label{sec:intro}
Consider the initial value problem,
\begin{equation}\label{eq:lin}
\left\{
\begin{aligned}
i\eps\d_t v^\eps +\frac{1}{2}\eps^2\Delta v^\eps &
=\frac{|x|^2}{2}v^\eps,\ \ \ \ (t,x)\in \R_+\times \R^n, \\
v^\eps_{\mid t=0}& = f(x),
\end{aligned}
\right.
\end{equation}
where $\eps \in ]0,1]$ is a parameter going to zero and $f$ is a
smooth function, say $f\in {\Sch}(\R)$. The
potential is the isotropic harmonic potential,
\begin{equation}\label{def:V}
V(x)\equiv \frac{|x|^2}{2}=\frac{1}{2}(x_1^2 + \dots + x_n^2).
\end{equation}
The case of anisotropic harmonic potentials is discussed in
Sect.~\ref{sec:anis}. 
Even though no (rapid)
oscillation is present in the initial data, the solution $v^\eps$
is rapidly oscillating (at frequency $1/\eps$) for any positive time,
and focuses at time $t=\frac{\pi}{2}$ (Section~\ref{sec:bkw}). This
can be seen by a stationary phase argument applied to the Mehler's
formula (see \cite{Feyn}),
\begin{equation}\label{eq:mehler}
v^\eps(t,x)=\frac{1}{(2i\pi \eps \sin
t)^{n/2}}\int_{\R^n}e^{\frac{i}{\eps \sin t}
\left(\frac{|x|^2+|y|^2}{2}\cos t -x\cdot y \right)}f(y)dy=:U^\eps(t)f(x).
\end{equation}
Perturbations of the harmonic potential by other potentials
(sub-quadratic perturbation, see \cite{Zelditch83}, \cite{Fujiwara}, 
\cite{KRY}, or
super-quadratic perturbation, see \cite{Yajima96}) have been studied, and
in particular the role of these perturbations on the singularities of
the fundamental solution of the Schr\"odinger equation. 

In physics, nonlinear perturbations are considered, for Bose-Einstein
condensation (see \cite{CCT}), where the harmonic potential is used for
its confining properties,
$$ i\hbar\d_t \psi^\hbar +\frac{1}{2}\hbar^2\Delta \psi^\hbar 
=\frac{|x|^2}{2}\psi^\hbar +Ng |\psi^\hbar|^2 \psi^\hbar,$$
where $N$ stands for the number of particles and $g$ is a coupling
constant (in $\hbar^2$). 

We study precisely the perturbation of (\ref{eq:lin}) with a nonlinear
term, 
\begin{equation}\label{eq:nl}
\left\{
\begin{aligned}
i\eps\d_t u^\eps +\frac{1}{2}\eps^2\Delta u^\eps &
=\frac{|x|^2}{2}u^\eps +\eps^{n\sigma} |u^\eps|^{2\sigma} u^\eps,\ \ \ \
(t,x)\in 
\R_+\times \R^n, \\ 
u^\eps_{\mid t=0} & = f(x) +r^\eps(x),
\end{aligned}
\right.
\end{equation}
with $\sigma >1/n$ if $n=1,2$, and $\frac{2}{n+2}<\sigma<\frac{2}{n-2}$
if $n\geq 3$. We assume that the perturbation
$r^\eps$ of the initial data is small in 
\begin{equation}\label{def:Sigma}
\Sigma := H^1(\R^n)\cap {\F}\left(H^1(\R^n)\right),
\end{equation}
where the Fourier transform is defined by 
$${\F}v(\xi)=\widehat v(\xi)=\int_{\R^n} e^{-ix\cdot \xi}v(x)dx,$$
and that $f\in \Sigma$. 
The space $\Sigma$ is equipped with the norm
$$\|f\|_\Sigma =\|f\|_{L^2} + \|\nabla_x f\|_{L^2} +\|xf\|_{L^2},$$
and we assume $\|r^\eps\|_\Sigma \Tend \eps 0 0$. 

\rq
 Initial data with plane oscillations. Let $\xi_0\in \R^n$, and
introduce
$${\tt u}^\eps(t,x)=u^\eps(t,x-\xi_0 \sin t)e^{i\left( x-
\frac{\xi_0}{2}\sin t\right)\cdot\xi_0 \cos t /\eps}.$$
Then ${\tt u}^\eps$ solves the Schr\"odinger equation
(\ref{eq:nl}), with initial data
$${\tt u^\eps}_{\mid t=0}= \left(f(x) + r^\eps(x) 
\right)e^{i\frac{x\cdot\xi_0}{\eps}}.$$
Therefore, describing the solution of (\ref{eq:nl}) is
enough to describe the solution when the initial data have plane
oscillations. \\

We can also prove some results with a \emph{focusing}
critical nonlinearity ($2\sigma=4/n$),
\begin{equation}\label{eq:nlfoc}
\left\{
\begin{aligned}
i\eps\d_t u^\eps +\frac{1}{2}\eps^2\Delta u^\eps &
=\frac{|x|^2}{2}u^\eps -\eps^2 |u^\eps|^{4/n} u^\eps,\ \ \ \
(t,x)\in 
\R_+\times \R^n, \\ 
u^\eps_{\mid t=0} & = f(x).
\end{aligned}
\right.
\end{equation}
We will consider the focusing case only in the one-dimensional
situation, and state the corresponding results at the end of this
introduction. Similar results for the multi-dimensional case would be
easy to prove. 

The idea of this paper is the following. Initially, the nonlinear term
is negligible, essentially because the term $|u^\eps|^{2\sigma}$ is
uniformly bounded in suitable Lebesgue spaces, therefore it
vanishes in 
the limit $\eps\rightarrow 0$ because of the factor
$\eps^{n\sigma}$. Meanwhile, the harmonic potential makes the solution
focus near the origin at time $t=\pi/2$, as in the linear case
(\ref{eq:lin}). When the focusing effects become relevant, that is when
$u^\eps$ becomes of order $\eps^{-n/2}$, the nonlinear term is no
longer negligible. On the other hand, if $u^\eps$ is localized near
$x=0$, the term $x^2u^\eps$ becomes
negligible; only the nonlinear term is relevant near the focus. When
the nonlinearity is \emph{defocusing} (Eq.~(\ref{eq:nl})), the
solution $u^\eps$ passes through the focus, and the crossing is given by 
the (nonlinear) scattering operator
associated to the unscaled Schr\"odinger equation,  
\begin{equation}\label{eq:nls}
i\d_t \psi +\frac{1}{2}\Delta \psi =|\psi|^{2\sigma} \psi.
\end{equation}
Since the nonlinearity is defocusing, dispersive effects in
(\ref{eq:nls}) are the same as for the free equation. Therefore, the
solution $u^\eps$ leaves the focus along dispersive rays. When rays
are dispersed, the energy is no longer localized, the nonlinear term
becomes negligible again and the harmonic potential makes the rule, as
before the focus (Th.~\ref{th:princ}). 
This process can be iterated indefinitely, and each
focus crossing is described by the scattering operator
(Cor.~\ref{cor:iter}).   

When the nonlinearity is \emph{focusing} (Eq.~(\ref{eq:nlfoc})), and
when the mass of $f$ is critical (see \cite{MerleDuke}),
the
solution blows up near $t=\pi/2$ (before or after, see
Prop.~\ref{prop:blowup}). The focusing effects of the
harmonic potential first, then of the nonlinear term, cumulate and
ruin the existence of the solution (Prop.~\ref{prop:blowup}). 

In both situations (focusing or defocusing nonlinearity), two distinct
r\'egimes occur. First, the harmonic potential leads the evolution of
the solution, next the nonlinear term does so. The two dynamics
superpose: they balance each other in the case of a defocusing
nonlinearity, and cumulate in the case of a focusing nonlinearity. The
matching of these two r\'egimes occurs in a boundary layer of size
$\eps$ around the focus, as in \cite{BG3}, \cite{GG}. 

Formal WKB expansions suggest that with
our choice $n \sigma >1$, the nonlinear term is negligible so long as no
focusing occurs. We prove that this holds true. It would not be so
with the choice $n\sigma =1$; the nonlinear term would be nowhere
negligible, and we leave out this case.

On the other hand, we
show that  the nonlinear term 
alters the asymptotics of the exact 
solution near and past the (first) focus. More precisely, we prove that
the caustic crossing is measured by the scattering operator
associated to (\ref{eq:nls}).
This phenomenon is to be compared with the results of \cite{Ca2},
where focusing is caused by initial 
oscillations, and with the results of \cite{BG3} (see also \cite{BG},
\cite{BG2}), where such a behavior was first 
noticed, for the wave equation. In the present case, focusing is
caused by the  
oscillations created by the harmonic potential, but the description of
the phenomena near the focal point is similar.

The asymptotic state for (\ref{eq:nls}) we will consider is defined by
\begin{equation}\label{eq:psi-}
\psi_-(x):= \frac{1}{(2i\pi)^{n/2}}\widehat f (x).
\end{equation}
We assume that $f\in \Sigma$ and that the scattering operator $S$ acts on
$\psi_-$, with 
$\psi_+= S\psi_-\in \Sigma$ (see
Proposition~\ref{prop:scattnls}), which is verified in either of the
following cases,
\begin{itemize}
\item $\sigma > \frac{2-n +\sqrt{n^2+12n +4}}{4n}$, or
\item $\|f\|_\Sigma$ is sufficiently small.
\end{itemize}
\begin{hyp}\label{hyp}Our hypotheses are the following:
\begin{itemize}
\item $1\leq n\leq 5$ and $\sigma >1/2$, so that the nonlinearity
$|z|^{2\sigma}z $ is twice differentiable.
\item If $n=1$, we assume moreover that $\sigma >1$.
\item If $3\leq n\leq 5$, we take $\sigma<\frac{2}{n-2}$.
\item If $n\leq 2$, we assume \begin{itemize}
           \item either $\sigma>\dis \frac{2-n +\sqrt{n^2+12n +4}}{4n}$, 
            \item or $\|f\|_\Sigma\leq \delta$ sufficiently small.
           \end{itemize}
\end{itemize}
\end{hyp}
\rq We could treat the case $n\geq 6$ if we replaced the nonlinear
term $\eps^{n\sigma}|u^\eps|^{2\sigma} u^\eps$ by $F(\eps^n
|u^\eps|^2) u^\eps$, with $F$ smooth and 
$$ F(|z|^2) \lesssim 1 + |z|^{2\sigma}.$$

\begin{theo}\label{th:princ}
Let $2 < r< \frac{2n}{n-2}$. If $n=1$, take $r=\infty$. Then 
under Assumptions~\ref{hyp}, the following asymptotics holds in
$L^2\cap L^r$, 
\begin{itemize}
\item If $0\leq t<\pi/2$, then 
$$u^\eps(t,x) \Eq \eps 0
\frac{e^{in\frac{\pi}{4}}}{(2\pi |\cos 
t|)^{n/2}}\widehat{\psi_-}\left(\frac{-x}{\cos t} \right)
e^{-i\frac{|x|^2}{2\eps}\tan t}.$$
\item If $\pi/2< t<3\pi/2$, then 
$$u^\eps(t,x) \Eq \eps 0
\frac{e^{in\frac{\pi}{4}-in\frac{\pi}{2}}}{(2\pi |\cos
t|)^{n/2}}\widehat{\psi_+}\left(\frac{-x}{\cos t} \right)
e^{-i\frac{|x|^2}{2\eps}\tan t},$$
\end{itemize}
where $\psi_-$ is defined by (\ref{eq:psi-}) and $\psi_+=S\psi_-$. 
\end{theo}

\rq We will prove actually that these asymptotics hold in a stronger
sense (see Corollary~\ref{cor:idem}, Propositions~\ref{prop:avt} and
\ref{prop:apres}). \\

We can restate this result when time $t=\pi/2$ is considered as the
initial time, in place of $t=0$. 
\begin{cor}\label{cor:transl}Let $\varphi\in \Sigma$. 
Assume that $u^\eps$ solves
\begin{equation}\label{eq:tdiff}
\left\{
\begin{aligned}
i\eps\d_t u^\eps +\frac{1}{2}\eps^2\Delta u^\eps &
=\frac{|x|^2}{2}u^\eps +\eps^{n\sigma} |u^\eps|^{2\sigma} u^\eps,\ \ \ \
(t,x)\in 
\R\times \R^n, \\ 
u^\eps_{\mid t=0} & =
\frac{1}{\eps^{n/2}}\varphi\left(\frac{x}{\eps}\right)  +
\frac{1}{\eps^{n/2}}r^\eps\left(\frac{x}{\eps}\right),
\end{aligned}
\right.
\end{equation}
with $\|r^\eps\|_\Sigma \Tend \eps 0 0$, and $\varphi$ satisfies the
same assumptions as $f$.  
Denote $\psi_\pm = W^{-1}_\pm \varphi$, where $W_\pm$ are the
wave operators (see Proposition~\ref{prop:scattnls}). 
Then with $r$ as in Th.~\ref{th:princ} and under
Assumptions~\ref{hyp}, the following asymptotics 
holds in $L^2\cap L^r$,  
\begin{itemize}
\item If $0< t<\pi$, then 
$$u^\eps(t,x) \Eq \eps 0 \left(\frac{-i}{2\pi \sin
t}\right)^{n/2}\widehat{\psi_+}\left(\frac{x}{\sin t} \right)
e^{i\frac{|x|^2}{2\eps\tan t}}.$$
\item If $-\pi< t<0$, then 
$$u^\eps(t,x) \Eq \eps 0 \left(\frac{-i}{2\pi \sin
t}\right)^{n/2}\widehat{\psi_-}\left(\frac{x}{\sin t} \right)
e^{i\frac{|x|^2}{2\eps\tan t}}.$$
\end{itemize}
\end{cor}

\rq In \cite{NierENS}, the author considers equations which
can be compared to (\ref{eq:lin}), that is
\begin{equation}\label{eq:nier}
\left\{
\begin{aligned}
i\eps\d_t v^\eps +\frac{1}{2}\eps^2\Delta v^\eps &
=V(x)v^\eps +U\left( \frac{x}{\eps}\right)v^\eps,\\
v^\eps_{\mid t=0}& = \frac{1}{\eps^{n/2}}\varphi\left( \frac{x}{\eps}\right),
\end{aligned}
\right.
\end{equation}
where $U$ is a short range potential. The potential $V$ in that case
cannot be the harmonic potential, for it has to be bounded as well as
all its derivatives. In that paper, the author proved that under
suitable assumptions, the influence of $U$ occurs near $t=0$ and is
localized near the origin, while only the value $V(0)$ of $V$ at the origin
is relevant in this r\'egime. For times $\eps \ll |t|<T_*$, the
situation is different: the potential $U$ becomes negligible, while
$V$ dictates the propagation. As in our paper, the transition between
these two r\'egimes is measured by the scattering operator associated
to $U$. 

Our assumption $n \sigma >1$ makes the nonlinear term short range. With
our scaling for the nonlinearity, this perturbation is relevant only
near the focus, where the harmonic potential is negligible, while
the opposite occurs for $\eps\ll |t|<\pi$. In this perspective, a new
point in our paper (besides the fact that the problem is nonlinear) is
that we can tell what happens for \emph{any} time, as stated in the
following corollary.

\begin{cor}\label{cor:iter}
Suppose the Assumption~\ref{hyp} are satisfied. Let $k\in
\N^*$. Then, with $r$ as in Th.~\ref{th:princ},
the asymptotics of $u^\eps$ for $\pi/2 +
(k-1)\pi < t < \pi/2 +k\pi$ is given, in $L^2\cap
L^r$, by
$$u^\eps(t,x) \Eq \eps 0 
\frac{e^{in\frac{\pi}{4}-ink\frac{\pi}{2}}}{(2\pi |\cos
t|)^{n/2}}\widehat{S^k\psi_-}\left(\frac{-x}{\cos t} \right)
e^{-i\frac{|x|^2}{2\eps}\tan t},$$
where $S^k$ denotes the $k$-th iterate of $S$ (which is well defined
under our assumptions on $f$). 
\end{cor}

\rq The phase shift $e^{in\frac{\pi}{4}-ink\frac{\pi}{2}}$ is present in
the linear case, for Eq.~(\ref{eq:lin}), and is explained in
\cite{Du}. On the contrary, the presence of the scattering operator
$S$ is typically a nonlinear phenomenon, as in \cite{Ca2}. The new
point here is that this operator is iterated, at each focus
crossing.\\

\rq If the nonlinear perturbation was of the form $\eps^{n\sigma_1}
|u^\eps|^{2\sigma_2} u^\eps$, with $\sigma_1 > \sigma_2>0$ (no
additional assumption on $\sigma_2$) and $n\sigma_1 >1$, the nonlinear term
would be everywhere negligible, that is, $S^k$ should be replaced by
the identity in Corollary~\ref{cor:iter}. This can be seen by an easy
adaptation of the proof of Theorem~\ref{th:princ}. This shows that the
scaling (\ref{eq:nl}) is critical for the nonlinearity to have a leading
order influence near the singularities ($t=\pi/2+k\pi$). 

We conclude this introduction by stating our result when the
nonlinearity is focusing (Eq.~(\ref{eq:nlfoc})).
\begin{prop}\label{prop:blowup}
Let $n=1$ and let $R$ be the unique solution  (up to translation and sign
change) of
$-\frac{1}{2}R'' +R =  R^5,$
given by,
\begin{equation}
  R(x)=\frac{3^{1/4}}{\sqrt{\cosh \left(2x\sqrt 2\right)}}.
\end{equation}
For $t_*\in\R$, define $f$ by
$f(x)=R(x)e^{i\frac{t_*}{2}x^2},$
and the approximate solution by
\begin{equation}\label{eq:solappfoc}
\tilde v^\eps (t,x)=\frac{1}{\sqrt{\frac{\pi}{2}+\eps t_*
-t}}R\left( \frac{x}{\frac{\pi}{2}+\eps t_*
-t} \right)e^{i\frac{\eps}{\pi/2 +\eps t_* -t}}e^{i\frac{|x|^2}{2\eps
(t-\pi/2 -\eps t_*)}}.
\end{equation}
Let $u^\eps$ be the solution of (\ref{eq:nlfoc}). Then for any
$\lambda > 0$,
$$\limsup_{\eps\rightarrow 0} \sup_{\frac{\pi}{2}-\Lambda \eps \leq t
\leq \frac{\pi}{2}+\eps t_* -\lambda\eps}\left\| B^\eps (t)\left(
u^\eps(t) - \tilde v^\eps (t)\right)\right\|_{L^2}\Tend \Lambda
{+\infty} 0,$$
where $B^\eps(t)$ is either of the operators $\Id$, $\eps\d_x$ or
$x/\eps +i(t-\pi/2)\d_x$.
In particular,
\begin{equation*} 
\begin{aligned}
\liminf_{\eps\rightarrow 0} \sup_{0\leq t\leq
\pi/2+t_*\eps-\lambda\eps} \|\eps\d_x u^\eps(t)\|_{L^2} &\Tend \lambda
{0^+} +\infty,\\
\liminf_{\eps\rightarrow 0} \sup_{0\leq t\leq
\pi/2+t_*\eps-\lambda\eps} \|\sqrt\eps u^\eps(t)\|_{L^\infty} &\Tend \lambda
{0^+} +\infty.
\end{aligned}
\end{equation*}
\end{prop}
\rq The blow up occurs at $t=\frac{\pi}{2}+\eps t_*$, no matter the
sign of $t_*$. This means that $u^\eps$ can blow up before or after
the focus.

This paper is organized as follows. In Sect.~\ref{sec:bkw}, we study
the linear equation (\ref{eq:lin}) using WKB methods, and introduce
two operators ($J^\eps$ and $H^\eps$) whose role is crucial in the
nonlinear setting. In Sect.~\ref{sec:asym}, we analyze the nonlinear
equation (\ref{eq:nl}), and we prove
Th.~\ref{th:princ}. In Sect.~\ref{sec:blowup}, we prove
Prop.~\ref{prop:blowup}. Finally, Sect.~\ref{sec:anis} addresses the case of 
anisotropic harmonic potentials. 

Some of the results presented in this paper were announced in
\cite{CaX01}.

\medskip

\noindent \emph{Acknowledgment.} The idea of studying such equations as
(\ref{eq:nl}) arose from discussions with L.~Miller, after an
invitation by P.~G\'erard and N.~Mauser at the Erwin
Schr\"odinger Institut in Vienna. This work was
partially supported by the ACI grant ``\'Equations des ondes :
oscillations, dispersion et contr\^ole'' and the START project
``Nonlinear Schr\"odinger 
equations'' of N.~Mauser. The  preliminary version of the manuscript was
improved thanks to remarks made by F.~Castella and B.~Bid\'egaray. 
   
\section{WKB expansion for the linear equation}
\label{sec:bkw}

\begin{figure}[htbp]
\begin{center}
\input{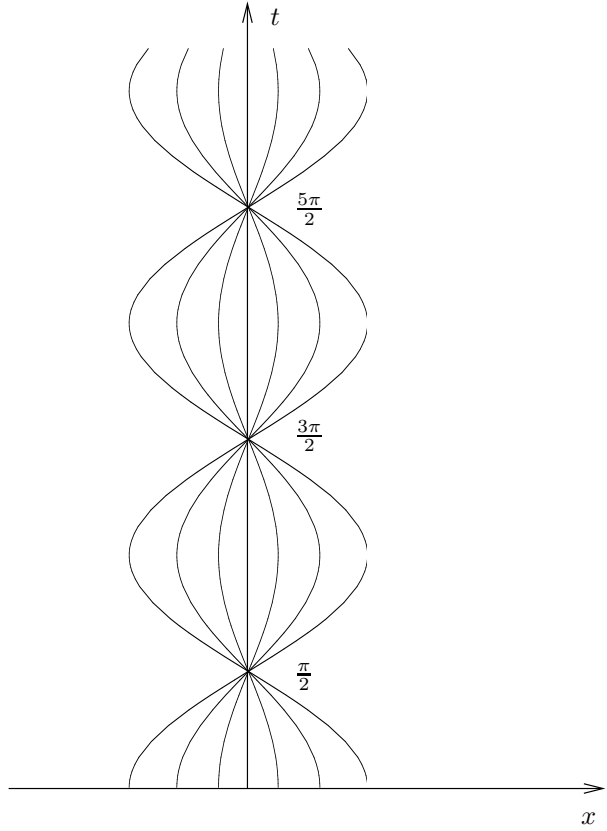}
\caption{Rays of geometric optics.}
\label{fig:sinus}
\end{center}
\end{figure}

We seek an approximate solution of the linear equation (\ref{eq:lin})
of the form,
\begin{equation}
v^\eps_{\rm app}(t,x)= v_0(t,x)e^{i\varphi(t,x)/\eps}.
\end{equation}
To cancel the term $\eps^0$ when plugging this
approximate solution in (\ref{eq:lin}), the phase $\varphi$ must
satisfy the eikonal equation,
\begin{equation}\label{eq:eikonal}
\d_t \varphi +\frac{1}{2}|\nabla_x \varphi|^2+\frac{|x|^2}{2}=0.
\end{equation}
To cancel the term $\eps^1$, the amplitude $v_0$ must satisfy the
transport equation, 
\begin{equation}\label{eq:transport}
\d_t v_0 +\nabla_x \varphi \cdot\nabla_x v_0 +\frac{1}{2}v_0 \Delta \varphi =0.
\end{equation}
To solve the eikonal equation, one computes the bicharacteristic
curves associated to the classical Hamiltonian
$$p(t,x,\tau,\xi)= \tau +\frac{|\xi|^2}{2}+\frac{|x|^2}{2},$$
given by
$$\left\{
\begin{array}{rcr}
\dot{t}&=&1,\\
\dot{x}&=&\xi,\\
\dot{\tau}&=&0,\\
\dot{\xi}&=&-x.
\end{array}
\right.$$
Therefore,
$$x(t)=x_0 \cos t +\xi_0 \sin t,\ \ \ \ \xi(t)=\xi_0 \cos t -x_0\sin
t.$$
Since no oscillation is present in the initial data, $\xi_0=0$, and
the rays of geometric optics are sinusoids,
\begin{equation}\label{eq:rays}
x(t)=x_0 \cos t.
\end{equation}
They all meet at the origin at time $t=\pi/2$, and periodically at time 
$t=\pi/2 +k\pi$ for any $k\in \N^*$ (Fig.~\ref{fig:sinus}).
 
Given $\xi(t) = \nabla_x \varphi (t)$, one can solve (\ref{eq:eikonal})
for $0\leq t< \pi /2$, by 
$$\varphi(t,x)=-\frac{|x|^2}{2}\tan t,$$
and (\ref{eq:transport}) is solved by 
$$v_0(t,x)=\frac{1}{(\cos t)^{n/2}}f\left(\frac{x}{\cos t} \right),$$
therefore
\begin{equation}\label{eq:approxlin}
v^\eps_{\rm app}(t,x)=\frac{1}{(\cos t)^{n/2}}f\left(\frac{x}{\cos t}
\right) e^{-i\frac{|x|^2}{2\eps}\tan t}. 
\end{equation}
Recall that $V(x)=\frac{|x|^2}{2}$. The approximate solution solves 
\begin{equation}\label{eq:linapp}
\left\{
\begin{aligned}
i\eps\d_t v^\eps_{\rm app} +\frac{1}{2}\eps^2\Delta v^\eps_{\rm app} &
=V(x)v^\eps_{\rm app}+\frac{1}{2}\eps^2e^{i\varphi(t,x)/\eps}\Delta v_0 ,\\
v^\eps_{{\rm{app}}\mid t=0}& = f(x).
\end{aligned}\right.
\end{equation}
Denote the remainder $w^\eps:=v^\eps-v^\eps_{\rm app}$. It solves,
\begin{equation}\label{eq:restelin}
\left\{
\begin{aligned}
i\eps\d_t w^\eps +\frac{1}{2}\eps^2\Delta w^\eps &
=V(x)w^\eps-\frac{1}{2}\eps^2e^{i\varphi(t,x)/\eps}\Delta v_0 ,\\
w^\eps_{\mid t=0}& =0.
\end{aligned}\right.
\end{equation}
Recall that 
$$\Delta v_0 (t,x)= \frac{1}{\cos^2 t}\frac{1}{(\cos t)^{n/2}}\Delta f\left
( \frac{x}{\cos t}\right).$$
Recall the classical result,
\begin{lem}\label{lem:estimL2}
Assume a function $w^\eps$ satisfies
\begin{equation}\label{eq:estimL2}
i\eps\d_t w^\eps +\frac{1}{2}\eps^2\Delta w^\eps 
=U(t,x)w^\eps + S^\eps(t,x), \ \ \ (t,x)\in I\times\R^n,
\end{equation}
where $U$ is a real-valued potential, $I$ is an interval, and
$S^\eps\in C_t(I,L^2)$. 
Then the following estimate holds for $t\in I$,
\begin{equation*}
\eps \d_t \|w^\eps(t)\|_{L^2} \leq 2 \|S^\eps(t)\|_{L^2}.
\end{equation*}
\end{lem}
Applying this lemma, it follows,
\begin{equation}\label{eq:estiml2}
\eps \d_t \|w^\eps(t)\|_{L^2}\leq \eps^2\frac{1}{\cos^2 t}\|\Delta f\|_{L^2}.
\end{equation}

With the idea of a nonlinear perturbation in mind, it is natural to seek
estimates in other spaces than $L^2$, in particular Sobolev like
spaces. In geometrical optics, it is classical to assess
$\eps$-derivatives to get nonlinear estimates (see for instance
\cite{RauchUtah}). This is because $\eps$-oscillating solutions are
studied. This approach is sharp for multi-phase problems, but it
contains no geometric information (given by the phase(s)). In our
case, only one phase is present, and in the nonlinear setting
(\ref{eq:nl}), it remains so. In the linear case, this means that
controlling $v^\eps_{\rm app}$ in Lebesgue's spaces $L^p$ is
equivalent to controlling $v_0$ in  $L^p$. With Gagliardo-Nirenberg
inequalities in mind, it is therefore natural to introduce the
operator
\begin{equation}
\boxed{J^\eps(t)= -i(\cos t) e^{i\varphi/\eps}\nabla_x
(e^{-i\varphi/\eps}\cdot)=  
\frac{x}{\eps}\sin t -i \cos t \nabla_x.} 
\end{equation}
Given the dynamics of the harmonic potential, it is also natural to
introduce the ``orthogonal'' operator,
\begin{equation}\boxed{
H^\eps(t)= x\cos t+i\eps \sin t \nabla_x.}
\end{equation}
When $n\geq 2$, we write, for $1\leq j\leq n$, 
\begin{equation}\begin{aligned}
H_j^\eps(t)&= x_j\cos t+i\eps \sin t \d_{x_j},\\
J_j^\eps(t)&= \frac{x_j}{\eps}\sin t -i \cos t \d_{x_j}.
                \end{aligned}
\end{equation}
We now state all the properties we will need, including the action on
nonlinear terms. 
\begin{lem}\label{lem:operators}
The operators $H^\eps$ and  $J^\eps$ satisfy the following properties.
\begin{itemize}
\item The commutation relation,
\begin{equation}\label{eq:commut}
\left[H_j^\eps(t), i\eps \d_t +\frac{1}{2}\eps^2\Delta -\frac{|x|^2}{2}
\right]=\left[J_j^\eps(t), i\eps \d_t +\frac{1}{2}\eps^2\Delta -\frac{|x|^2}{2}
\right]=0.
\end{equation}
\item Denote $M^\eps(t)= e^{-i\frac{|x|^2}{2\eps}\tan t}$, and 
$Q^\eps(t) = e^{i\frac{|x|^2}{2\eps\tan t}}$,
then
$J^\eps(t)$  and $H^\eps(t)$ read, for $t\not \in \frac{\pi}{2}\Z$,
\begin{equation}\label{eq:factor}
J^\eps(t) =-i\cos t M^\eps(t)\nabla_x  M^\eps(-t), \ \ 
H^\eps(t) =i\eps\sin t Q^\eps(t)\nabla_x  Q^\eps(-t)
\end{equation}
\item The modified Sobolev inequalities. For $n=1$ and $t\not \in
\frac{\pi}{2}\Z$, 
\begin{equation}\label{eq:sobolev}
\begin{aligned}
\|w(t)\|_{L^\infty} & \leq \frac{C}{|\cos
t|^{1/2}}\|w(t)\|_{L^2}^{1/2}
\|J^\eps(t)w(t)\|_{L^2}^{1/2}, \\
\|w(t)\|_{L^\infty} & \leq \frac{C}{|\eps \sin
t|^{1/2}}\|w(t)\|_{L^2}^{1/2}
\|H^\eps(t)w(t)\|_{L^2}^{1/2}.
\end{aligned}
\end{equation}
For $n\geq 2$, and $2\leq r< \frac{2n}{n-2}$, define $\delta(r)$ by
$$\delta(r)\equiv n\left(\frac{1}{2}-\frac{1}{r}\right).$$
Then for any $2\leq r<\frac{2n}{n-2}$, there exists $C_r$ such
that, for $t\not \in \frac{\pi}{2}\Z$,
\begin{equation}\label{eq:sobolevn}
\begin{aligned}
\|w(t)\|_{L^r} & \leq \frac{C_r}{|\cos
t|^{\delta(r)}}\|w(t)\|_{L^2}^{1-\delta(r)}
\|J^\eps(t)w(t)\|_{L^2}^{\delta(r)}, \\
\|w(t)\|_{L^r} & \leq \frac{C_r}{|\eps \sin
t|^{\delta(r)}}\|w(t)\|_{L^2}^{1-\delta(r)}
\|H^\eps(t)w(t)\|_{L^2}^{\delta(r)}.
\end{aligned}
\end{equation}
\item For any function $F\in C^1(\C, \C)$ satisfying the gauge
invariance condition 
$$\exists G\in C(\R_+,\R),\ F(z)=zG(|z|^2),$$
one has, for $t\not \in \frac{\pi}{2}\Z$,
\begin{equation}\label{eq:jauge}\begin{aligned}
H^\eps(t)F(w)&= \d_zF(w)H^\eps(t)w - \d_{\bar
z}F(w)\overline{H^\eps(t)w},\\
J^\eps(t)F(w)&= \d_zF(w)J^\eps(t)w - \d_{\bar
z}F(w)\overline{J^\eps(t)w}.
                                \end{aligned}
\end{equation} 
\end{itemize}
\end{lem}
\rqs
\begin{itemize}
\item Estimates (\ref{eq:sobolev}) are easy consequences of the
conjugation properties (\ref{eq:factor}).
\item With the WKB approximation (\ref{eq:approxlin}) in mind, the $|\cos
t|^{-1/2}$ term in (\ref{eq:sobolev}) gives optimal time dependence of
the $L^\infty_x$ estimates of the solution of (\ref{eq:lin}) away from
the focus. This is the main advantage of this operator over all the
others one could think of (such as $\nabla_x$ in particular). 
\item The $|\eps \sin t|^{-1/2}$ term in (\ref{eq:sobolev}) gives
optimal $L^\infty_x$ estimates of the solution of (\ref{eq:lin}) near
the focus (where $|\sin t|\sim 1$). 
\item The
operator $J^\eps$ can be considered as the modification of the
Galilean operator $x+it\nabla_x$, which is very useful in scattering
theory (see 
\cite{Caz}, 
\cite{GinibreDEA}, \cite{Ginibre}).
For semi-classical problems where focusing at the origin occurs, it
was used in \cite{Ca2} and
\cite{CaCMP}, with the rescaling $\frac{x}{\eps}+i(t-t_*)\nabla_x$, where
$t_*$ is the focusing time. The operator $J^\eps$ is that operator,
transported to the case of a harmonic potential. 
\item Property (\ref{eq:jauge}) states that $H^\eps$ and $J^\eps$ act
on nonlinearities satisfying the gauge
invariance condition like derivatives (Eq.~(\ref{eq:jauge}) holds for
the operator $i\nabla_x$).
\item The fact that all these identities, except the first one,  hold
only for almost all $t\in \R$ is not a problem, since in any case
integrations with respect to time will be performed.
\item The operators $J^\eps$ and $H^\eps$ are known in quantum
mechanics, as Heisenberg observables (metaplectic transforms, see
 \cite{GS}, \cite{Folland}),
\begin{eqnarray}
J^\eps(t) & = & U^\eps(t)(-i\nabla_x) U^\eps(-t) \label{eq:num1} \\
          & = & U^\eps\left(t-\frac{\pi}{2}\right)\frac{x}{\eps}
          U^\eps\left(\frac{\pi}{2}-t\right) \label{eq:num2},\\
H^\eps(t) & = & U^\eps(t) x U^\eps(-t) \label{eq:num3} \\
          & = & U^\eps\left(t-\frac{\pi}{2}\right)(i\eps \nabla_x)
          U^\eps\left(\frac{\pi}{2}-t\right) \label{eq:num4},
\end{eqnarray}
where $U^\eps(t)$ is the propagator defined by Mehler's formula 
(\ref{eq:mehler}), that is
$$U^\eps(t)= e^{-i\frac{t}{2\eps}\left( -\eps^2\Delta
+|x|^2\right)} .$$ 
The commutation properties (\ref{eq:commut}) are straightforward
consequences of the conjugation relations (\ref{eq:num1}) and
(\ref{eq:num3}). Identities between (\ref{eq:num1}) and
(\ref{eq:num2}) on the
one hand, (\ref{eq:num3}) and  (\ref{eq:num4}) on the other hand, 
are due to the geometric properties of the harmonic oscillator, that
rotates the phase space. It is
easy to check that Sobolev inequalities
 (\ref{eq:sobolev}) follow from
(\ref{eq:num2}), (\ref{eq:num3}) and the estimate
$$\|U^\eps(t)f\|_{L^\infty_x}\lesssim \frac{1}{(\eps |\sin
t|)^{1/2}}\|f\|_{L^1}.$$
The most remarkable fact is certainly that in the case of the harmonic
potential, one can estimate the action of these observables of Heisenberg
on a large class of nonlinearities, through (\ref{eq:jauge}).
\end{itemize}
Lemma \ref{lem:operators} makes it possible to get more precise
estimates of the approximation given by (the first term of) WKB
methods. Denote 
\begin{equation}\label{eq:defH}
{\EuScript H}:= \{f\in H^3(\R^n), \textrm{ such that } xf \in H^2(\R^n) \}.
\end{equation}
\begin{prop}\label{prop:bkw}
Assume $f\in {\EuScript H}$. Then there exists
$C=C(\|f\|_{H^3}, \|xf\|_{H^2})$ such that 
the remainder $v^\eps-v^\eps_{\rm app}$ satisfies, for $0\leq
t<\pi/2$, 
\begin{equation*}
\begin{aligned}
\|(v^\eps-v^\eps_{\rm app})(t)\|_{L^2}+&\|J^\eps(v^\eps-v^\eps_{\rm
app})(t)\|_{L^2}+ \|H^\eps(v^\eps-v^\eps_{\rm app})(t)\|_{L^2} \leq \\
& \leq C \left(\eps 
\int_0^t \frac{ds}{\cos^2 s}+ \frac{\eps^2}{\cos^2 t} \right).
\end{aligned}
\end{equation*}
\end{prop}
\rq As mentioned in the introduction (Eq.~(\ref{eq:mehler}), the
expression of $v^\eps$ is  given explicitly 
by an oscillatory integral, and the above result could be proved by
stationary phase methods. Nevertheless, we do not use this approach,
and rather present the approach whose spirit is the same as in the
nonlinear setting. 

\begin{proof}
The first estimate is given by (\ref{eq:estiml2}), with
$C=\|\Delta f\|_{L^2}$. For the second estimate, apply $J^\eps(t)$ to
(\ref{eq:restelin}). The commutation property (\ref{eq:commut}) yields,
\begin{equation}\label{eq:restelin2}
\left\{
\begin{aligned}
i\eps\d_t J^\eps w^\eps +\frac{1}{2}\eps^2\Delta J^\eps w^\eps &
=V(x)J^\eps
w^\eps-\frac{1}{2}\eps^2J^\eps(t)\left(e^{i\varphi(t,x)/\eps}\Delta 
v_0  \right),\\
J^\eps w^\eps_{\mid t=0}& =0.
\end{aligned}\right.
\end{equation}
One has explicitly,
\begin{equation*}
\begin{aligned}
J^\eps(t)\left(e^{i\varphi(t,x)/\eps}\Delta 
v_0  \right)& = \left(\frac{x}{\eps}\sin t - i\cos t \nabla_x
\right)\left( e^{-i\frac{|x|^2}{2\eps}\tan t}\frac{1}{(\cos
t)^{n/2+2}}\Delta f\left(\frac{x}{\cos t} \right)\right)\\ 
&= - ie^{-i\frac{|x|^2}{2\eps}\tan t}
\cos t \nabla_x \left( \frac{1}{(\cos
t)^{n/2+2}}\Delta f\left(\frac{x}{\cos t} \right)\right)\\
&= - ie^{-i\frac{|x|^2}{2\eps}\tan t}\frac{1}{(\cos
t)^{n/2+2}}\nabla_x\Delta f\left(\frac{x}{\cos t} \right),
\end{aligned}
\end{equation*}
and the same estimate as for the $L^2$ case follows, with $C=
\|f\|_{H^3}$. For the last estimate of the proposition, apply $H^\eps(t)$ to
(\ref{eq:restelin}). Because of the commutation property
(\ref{eq:commut}), the remainder $H^\eps(t)w^\eps$ is estimated by the
$L^2$ norm of 
\begin{equation*}
\begin{aligned}
H^\eps(t)\left(e^{i\varphi(t,x)/\eps}\Delta
v_0  \right)=& \left(x\cos t + i\eps\sin t \nabla_x
\right)\left( e^{-i\frac{|x|^2}{2\eps}\tan t}\frac{1}{(\cos
t)^{n/2+2}}\Delta f\left(\frac{x}{\cos t} \right)\right)\\ 
 =& \frac{x}{\cos t} e^{-i\frac{|x|^2}{2\eps}\tan t}\frac{1}{(\cos
t)^{n/2+2}}\Delta f\left(\frac{x}{\cos t} \right) \\
& + i\eps \tan t \frac{1}{(\cos
t)^{n/2+2}}\nabla_x\Delta f\left(\frac{x}{\cos t} \right).
\end{aligned}
\end{equation*}
The $L^2$ norm of the first term is $\frac{1}{\cos^2
t}\|x\Delta f\|_{L^2}$,  and the $L^2$ norm of the second term is
$\eps\frac{\sin t}{\cos^3
t}\|f\|_{\dot{H}^3}$. This completes the proof of
Proposition~\ref{prop:bkw}. 
\end{proof}

\smallbreak

From Prop.~\ref{prop:bkw}, WKB methods provide a good
approximation of the exact solution before focusing. More precisely,
the remainder will be small up to a boundary layer of size $\eps$
around $t=\pi/2$. 

The assumption $f\in {\EuScript H}$ is necessary to
estimate precisely the validity of WKB approximation, but is not really
essential. Since the set of such $f$ is dense in $\Sigma$, the
following lemma shows that this extra regularity can be
introduced without modifying the asymptotics. 
\begin{lem}
Assume $f\in \Sigma$, and let $v^\eps$ be the solution of (\ref{eq:lin}). 
Then for any $t>0$,
$$\|v^\eps(t)\|_{L^2}=\|f\|_{L^2};\ \
\|J^\eps(t)v^\eps\|_{L^2}=\|\nabla f\|_{L^2};\ \
\|H^\eps(t)v^\eps\|_{L^2}=\|xf\|_{L^2}.$$
\end{lem}
\begin{proof}
This lemma is a straightforward consequence of
Lemma~\ref{lem:estimL2} and of the commutation property
(\ref{eq:commut}). 
\end{proof}

Notice that the $L^2$-norm of $v^\eps_{\rm app}(t)$ does not depend on
time, nor that of $J^\eps(t)v^\eps_{\rm app}$ or $H^\eps(t)v^\eps_{\rm
app}$. We can therefore remove the smoothness assumption of
Prop.~\ref{prop:bkw}.  
\begin{cor}\label{cor:idem}
Assume $f\in \Sigma$. Then,
\begin{equation*}
\limsup_{\eps \rightarrow 0} \sup_{0\leq t\leq\frac{\pi}{2}-\Lambda
\eps}\left\|A^\eps(t)\left( v^\eps-v^\eps_{\rm
app})(t)\right)\right\|_{L^2} \Tend \Lambda {+\infty} 0,
\end{equation*}
where $A^\eps(t)$ is either of the operators $\Id$, $J^\eps(t)$ or
$H^\eps(t)$.  
\end{cor}

\section{The nonlinear case}
\label{sec:asym}
The proof for asymptotics in the nonlinear setting relies on
Strichartz estimates (even though we could do without when $n=1$). We
first recall how we get them in the present case, then prove a general
estimate. Then the proof of Th.~\ref{th:princ} is essentially split
into three parts: the asymptotics before the focus ($0\ll \pi/2-t$),
the matching between the two r\'egimes (linear and nonlinear), and the
asymptotics around the focus ($|t-\pi/2|\lesssim \eps$). 
\subsection{Strichartz inequalities}
\label{sec:stri}
First, recall the classical definition (see e.g. \cite{Caz}),
\begin{defin}\label{def:adm}
 A pair $(q,r)$ is {\bf admissible} if $2\leq r
  <\frac{2n}{n-2}$ (resp. $2\leq r\leq \infty$ if $n=1$, $2\leq r<
  \infty$ if $n=2$)  
  and 
$$\frac{2}{q}=\delta(r)\equiv n\left( \frac{1}{2}-\frac{1}{r}\right).$$
\end{defin}
Strichartz estimates provide mixed type estimates (that is, in spaces
of the form $L^q_t(L^r_x)$ with $(q,r)$ admissible) of quantities
involving the unitary group $e^{i\frac{t}{2}\Delta}$ (see \cite{Strichartz},
\cite{GV85}, \cite{Kato87},  \cite{Yajima87}, \cite{Caz}, \cite{GinibreDEA},
\cite{Ginibre}). With the scaling  
of Eq.~(\ref{eq:nl}), the natural unitary group to consider is
\begin{equation}\label{eq:U0eps}
U_0^\eps(t):=e^{i\eps \frac{t}{2}\Delta}.
\end{equation}
Now we can state the Strichartz estimates obtained by a scaling
argument from the usual ones (with $\eps =1$). The notation $r'$
stands for the H\"older conjugate exponent of $r$.
\begin{prop}\label{prop:stri1}(Scaled Strichartz inequalities)
  \begin{enumerate}
  \item For any admissible pair $(q,r)$, there exists $C_r$ such that
\begin{equation}\label{eq:strichlib}
   \eps^{\frac{1}{q}} \left\| U^\eps_0(t)u\right\|_{L^q(\R;L^r)}\leq C_r
   \|u\|_{L^2}. 
  \end{equation}
  \item For any admissible pairs $(q_1,r_1)$ and $(q_2,r_2)$, and any
    interval $I$, there exists $C_{r_1,r_2}$ such that 
\begin{equation}\label{eq:strichnl}
      \eps^{\frac{1}{q_1}+\frac{1}{q_2}}\left\| \int_{I\cap\{s\leq
      t\}} U^\eps_0(t-s)F(s)ds 
      \right\|_{L^{q_1}(I;L^{r_1})}\leq C_{r_1,r_2} \left\|
      F\right\|_{L^{q'_2}(I;L^{r'_2})}. 
    \end{equation}
  \end{enumerate}
 The above constants are independent of $\eps$ and $I$. 
\end{prop}
The proof of this result relies on two properties (see \cite{Caz}, or
\cite{KT} for a more general statement):
\begin{itemize}
\item The group $U_0^\eps$ is unitary on $L^2$,
$\|U_0^\eps(t)\|_{L^2\rightarrow L^2} =1$.
\item For $t\not = 0$, it maps $L^1(\R^n)$ into $L^\infty(\R^n)$,
$$\|U_0^\eps(t)\|_{L^1\rightarrow L^\infty}\lesssim \frac{1}{(\eps
|t|)^{n/2}}.$$ 
\end{itemize}
As a matter of fact, these two estimates also hold for the propagator
associated to the Schr\"odinger equation with a harmonic potential
(\ref{eq:lin}). Therefore we can obtain similar Strichartz estimates
(see \cite{Caz}). 

If $v^\eps$ solves (\ref{eq:lin}), then Mehler's formula yields, for $t\not
\in \pi\Z$ (see \cite{Feyn}), 
\begin{equation*}
v^\eps(t,x)=\frac{1}{(2i\pi \eps \sin
t)^{n/2}}\int_{\R^n}e^{\frac{i}{\eps \sin t}
\left(\frac{|x|^2+|y|^2}{2}\cos t -x\cdot y \right)}f(y)dy=:U^\eps(t)f(x).
\end{equation*}
Therefore:
\begin{itemize}
\item The group $U^\eps$ is unitary on $L^2$,
$\|U^\eps (t)f\|_{L^2} =\|f\|_{L^2}$.
\item For $t\in ]-\pi,0[\cup ]0,\pi[$, it maps $L^1(\R^n)$ into
$L^\infty(\R^n)$, 
$$\|U^\eps\|_{L^1\rightarrow L^\infty}\lesssim \frac{1}{(\eps
|\sin t|)^{n/2}}.$$ 
\end{itemize}
Since for $|t|\leq \pi/2$, $|\sin t|\geq \frac{2}{\pi}|t|$, the proof
of Prop.~\ref{prop:stri1} still works when $U_0^\eps$ is replaced by
$U^\eps$, provided that only \emph{finite} time intervals are
considered.
\begin{prop}\label{prop:stri2}
  \begin{enumerate}
  \item For any admissible pair $(q,r)$, for any finite interval $I$,
  there exists $C_r(I)$ such that 
\begin{equation}\label{eq:strichlib2}
   \eps^{\frac{1}{q}} \left\| U^\eps(t)u\right\|_{L^q(I;L^r)}\leq C_r(I)
   \|u\|_{L^2}. 
  \end{equation}
  \item For any admissible pairs $(q_1,r_1)$ and $(q_2,r_2)$, and any
   finite interval $I$, there exists $C_{r_1,r_2}(I)$ such that 
\begin{equation}\label{eq:strichnl2}
      \eps^{\frac{1}{q_1}+\frac{1}{q_2}}\left\| \int_{I\cap\{s\leq
      t\}} U^\eps_0(t-s)F(s)ds 
      \right\|_{L^{q_1}(I;L^{r_1})}\leq C_{r_1,r_2}(I) \left\|
      F\right\|_{L^{q'_2}(I;L^{r'_2})}. 
    \end{equation}
  \end{enumerate}
 The above constants are independent of $\eps$. 
\end{prop}
\subsection{A general estimate}
\label{sec:gen}
We start with an algebraic lemma.
\begin{lem}\label{lem:alg}Let $n\geq 2$, and
assume $\frac{2}{n+2}<\sigma<\frac{2}{n-2}$. 
There exists $\underline{q}$, $\underline{r}$, $\underline{s}$ and
$\underline{k}$ satisfying
\begin{equation}\label{eq:holder}
  \left\{
  \begin{aligned}
    \frac{1}{\underline{r}'}&=\frac{1}{\underline{r}}+
\frac{2\sigma}{\underline{s}},\\ 
    \frac{1}{\underline{q}'}&=\frac{1}{\underline{q}}+
\frac{2\sigma}{\underline{k}},
  \end{aligned}\right.
\end{equation}
and the additional conditions:
\begin{itemize}
\item The pair $(\underline{q},\underline{r})$ is admissible,
\item $0<\frac{1}{\underline{k}}<\delta(\underline{s})<1$. 
\end{itemize}
If $n=1$, we take $(\underline{q},\underline{r})=(\infty, 2)$ and
$(\underline{k}, \underline{s})=(2\sigma,\infty)$. 
\end{lem}
\begin{proof}
With
$\delta(\underline{s})=1$, the first part of  (\ref{eq:holder})
becomes 
$$\delta(\underline{r})=\sigma \left( \frac{n}{2}-1\right),$$
and this expression is less than $1$ for 
$\sigma<\frac{2}{n-2}$. Still with $\delta(\underline{s})=1$, the
second part of (\ref{eq:holder}) yields
$$\frac{2}{\underline{k}}=1-\frac{n}{2}+\frac{1}{\sigma},$$
which lies in $]0,2[$ for
$\frac{2}{n+2}<\sigma<\frac{2}{n-2}$. By continuity, these conditions
are still satisfied for $\delta (\underline{s})$ close to $1$ and $\delta
(\underline{s})<1$. 
\end{proof}

From now on, we assume $n\geq 2$ and
$\frac{2}{n+2}<\sigma<\frac{2}{n-2}$.
We state a general estimate that can be applied to nonlinear
Schr\"odinger equations with or without harmonic potential. Let ${\EuScript
U}^\eps(t)$ be a group for which Prop.~\ref{prop:stri2} holds
(typically, $U_0^\eps$ or $U^\eps$ in our situation). We seek a
general estimate for the integral equation,
\begin{equation}\label{eq:Duhgen}
\begin{aligned}
u^\eps(t)={\EuScript U}^\eps(t-t_0)u_0^\eps &-i \eps^{n\sigma -1}
\int_{t_0}^t {\EuScript U}^\eps(t-s)F^\eps (u^\eps)(s) ds \\
&-i\eps^{-1}
 \int_{t_0}^t {\EuScript U}^\eps(t-s)h^\eps (s) ds.
\end{aligned}
\end{equation}
This equation generalizes the Duhamel formula for Eq.~(\ref{eq:nl}),
\begin{itemize}
\item to the case of the same equation without potential (take
$U_0^\eps$ in place of $U^\eps$),
\item to the case of any initial time and any initial data ($u_0^\eps$
and $t_0$ are general), 
\item to the possibility of having a nonlinear term which is not a
power, $F^\eps(u^\eps)$,
\item to the possibility of having a source term, $h^\eps$. 
\end{itemize}
\begin{prop}\label{prop:estgen}Let $t_1>t_0$, with $|t_1-t_0|\leq \pi$. 
Assume that there exists a constant $C$ independent of $t$ and $\eps$
such that for $t_0\leq t\leq t_1$,
\begin{equation}\label{eq:hypF}
\|F^\eps(u^\eps)(t)\|_{L^{{\underline r}'}_x}\leq \frac{C}{\left(|\cos
t|+\eps \right)^{2\sigma
\delta(\underline{s})}}\|u^\eps(t)\|_{L^{\underline r}_x}\ \ \ ,
\end{equation}
and define
$$A^\eps(t_0,t_1):= \left(\int_{t_0}^{t_1} 
\frac{dt}{\left(|\cos t|+\eps \right)^{\underline{k}
\delta(\underline{s})}}\right)^{2\sigma/\underline{k}}.$$
Then there exist $C^*$ independent of $\eps$, $t_0$ and $t_1$ such
that for any admissible pair $(q,r)$, 
\begin{equation}\label{eq:estgen}
\begin{aligned}
\|u^\eps\|_{L^{\underline q}(t_0,t_1;L^{\underline r})} \leq & C^*
\eps^{-1/{\underline q}}\|u_0^\eps\|_{L^2} + C_{{\underline q}, q}
\eps^{-1-\frac{1}{{\underline q}}- \frac{1}{q}}\|h^\eps\|_
{ L^{q'}(t_0,t_1;L^{r'})}\\
& + C^*\eps^{2\sigma\left(\delta(\underline s)-\frac{1}{\underline
k}\right)} A^\eps(t_0,t_1) \|u^\eps\|_{L^{\underline
q}(t_0,t_1;L^{\underline r})}.
\end{aligned}
\end{equation}
\end{prop}
We will rather use the following corollary, 
\begin{cor}\label{cor:estgen}
Suppose the assumptions of Prop.~\ref{prop:estgen} are
satisfied. Assume moreover that
$C^*\eps^{2\sigma\left(\delta(\underline s)-\frac{1}{\underline 
k}\right)} A^\eps(t_0,t_1)\leq 1/2$, which holds
in either of the two cases, 
\begin{itemize}
\item $0\leq t_0\leq t_1\leq \frac{\pi}{2}-\Lambda \eps$, with
$\Lambda\geq \Lambda_0$ sufficiently large,
\item $t_0,t_1\in[\frac{\pi}{2}-\Lambda \eps,\frac{\pi}{2}+\Lambda
\eps]$, with $\frac{t_1-t_0}{\eps}\leq \eta$ sufficiently small. 
\end{itemize}
Then 
\begin{equation}\label{eq:estgenL2}
\|u^\eps\|_{L^\infty(t_0,t_1;L^2)} \leq  C
\|u_0^\eps\|_{L^2} + C_{{\underline q}, q}
\eps^{-1- \frac{1}{q}}\|h^\eps\|_
{ L^{q'}(t_0,t_1;L^{r'})}.
\end{equation}
\end{cor}

\noindent {\it Proof of Proposition~\ref{prop:estgen}.} Apply
Strichartz inequalities (\ref{eq:strichlib2}) and 
(\ref{eq:strichnl2}) to (\ref{eq:Duhgen}) with $q_1=\underline q$,
$r_1=\underline r$, and $q_2=\underline q$, $r_2=\underline r$
for the term with $F^\eps(u^\eps)$, $q_2=q$, $r_2=r$
for the term with $h^\eps$,  it yields 
\begin{equation*}
\begin{aligned}
\|u^\eps\|_{L^{\underline q}(t_0,t_1;L^{\underline r})} \leq & C
\eps^{-1/{\underline q}}\|u_0^\eps\|_{L^2} + C_{{\underline q}, q}
\eps^{-1-\frac{1}{{\underline q}}- \frac{1}{q}}\|h^\eps\|_
{ L^{q'}(t_0,t_1;L^{r'})}\\
& + C\eps^{n\sigma -1 -\frac{2}{\underline
q}}\|F^\eps(u^\eps)\|_{L^{{\underline q}'}(t_0,t_1;L^{{\underline r}'})}. 
\end{aligned}
\end{equation*}
Then estimate the space norm of the last term by (\ref{eq:hypF}) and
apply H\"older inequality in time, thanks to  (\ref{eq:holder}), it yields
(\ref{eq:estgen}). \qed  

\noindent {\it Proof of Corollary~\ref{cor:estgen}.} The additional
assumption implies that the last term in (\ref{eq:estgen}) can be
``absorbed'' by the left-hand side, up to doubling the constants, 
\begin{equation}\label{eq:abs!}
\|u^\eps\|_{L^{\underline q}(t_0,t_1;L^{\underline r})} \leq  C
\eps^{-1/{\underline q}}\|u_0^\eps\|_{L^2} + C
\eps^{-1-\frac{1}{{\underline q}}- \frac{1}{q}}\|h^\eps\|_
{ L^{q'}(t_0,t_1;L^{r'})}.
\end{equation}
Now apply Strichartz inequalities (\ref{eq:strichlib2}) and 
(\ref{eq:strichnl2}) to (\ref{eq:Duhgen}) again, but with $q_1=\infty$,
$r_1=2$, and $q_2=\underline q$, $r_2=\underline r$
for the term with $F^\eps(u^\eps)$, $q_2=q$, $r_2=r$
for the term with $h^\eps$. It yields 
\begin{equation*}
\begin{aligned}
\|u^\eps\|_{L^\infty(t_0,t_1;L^2)} \leq & C
\|u_0^\eps\|_{L^2} + C
\eps^{-1- \frac{1}{q}}\|h^\eps\|_
{ L^{q'}(t_0,t_1;L^{r'})}\\
& + C\eps^{n\sigma -1 -\frac{1}{\underline
q}}\|F^\eps(u^\eps)\|_{L^{{\underline q}'}(t_0,t_1;L^{{\underline r}'})}. 
\end{aligned}
\end{equation*}
Like before,
\begin{equation*}
\begin{aligned}
\eps^{n\sigma -1 -\frac{1}{\underline
q}}\|F^\eps(u^\eps)\|_{L^{{\underline q}'}(t_0,t_1;L^{{\underline
r}'})}&  \leq C \eps^{\frac{1}{\underline
q}}\eps^{2\sigma\left(\delta(\underline s)-\frac{1}{\underline
k}\right)} A^\eps(t_0,t_1) \|u^\eps\|_{L^{\underline
q}(t_0,t_1;L^{\underline r})}\\
& \leq C \eps^{\frac{1}{\underline
q}}\|u^\eps\|_{L^{\underline
q}(t_0,t_1;L^{\underline r})},
\end{aligned}
\end{equation*}
and the corollary follows from (\ref{eq:abs!}). \qed
\subsection{Existence results}
\label{sec:ex}
Local existence in $\Sigma$ stems from
the well-known case of the nonlinear Schr\"odinger equation
(\ref{eq:nls}), once we noticed that the operators $H^\eps$ and
$J^\eps$ are the exact substitutes for the usual operators $\eps\nabla$
and $\frac{x}{\eps}+i(t-\frac{\pi}{2})\nabla$, by
Lemma~\ref{lem:operators}. Duhamel's formula for (\ref{eq:nl}) writes
\begin{equation}\label{eq:Duh}
u^\eps(t)=U^\eps(t)(f+r^\eps) -i \eps^{n\sigma -1}
\int_0^t U^\eps(t-s)( |u^\eps|^{2\sigma}u^\eps)(s) ds.
\end{equation}
Replacing $U^\eps$ with $U^\eps_0$ would yield the Duhamel's formula
for the same equation with no harmonic potential. From the above
remark (the essential point is that $H^\eps$ and
$J^\eps$ commute with $U^\eps$)
and the fact that the same Strichartz inequalities hold for
$U^\eps$ and $U^\eps_0$ when time is bounded, local existence is
actually a byproduct of the existence theory for (\ref{eq:nl}) (which
relies essentially on the results of Sect.~\ref{sec:gen}, see
\cite{Kato87}, \cite{Caz}, \cite{GinibreDEA}, \cite{Ginibre}). For $(q_0,r_0)$
admissible, introduce the spaces 
\begin{equation*}
\begin{aligned}
Y^\eps_{r_0} (I)& =\{ u^\eps \in C(I,\Sigma), u^\eps, H^\eps
u^\eps, J^\eps u^\eps \in  
L^{q_0}_{loc}(I,L^{r_0}_x) \},\\
Y^\eps (I) &=\{ u^\eps \in C(I,\Sigma), \forall
(q,r)\textrm{ admissible}, u^\eps, H^\eps u^\eps, J^\eps u^\eps \in 
L^q_{loc}(I,L^r_x) \}. 
\end{aligned}
\end{equation*}

\begin{prop}\label{prop:exist}
Fix $\eps \in ]0,1]$, and let $f,r^\eps\in \Sigma$. There exists
$t^\eps >0$ such that (\ref{eq:nl}) has a unique solution $u^\eps \in
Y^\eps_{2\sigma+2}(0,t^\eps)$. Moreover, this solution belongs to
$Y^\eps(0,t^\eps)$.  The same result holds for Eq.~(\ref{eq:nlfoc})
and for any initial time.  
\end{prop}
We can take $t^\eps =+\infty$ when the
nonlinearity is defocusing (Eq.~(\ref{eq:nl})), thanks to the
conservations of mass and energy,
\begin{equation}\label{eq:masse}
\|u^\eps(t)\|_{L^2}=\|u^\eps(0)\|_{L^2}=O(1),
\end{equation}
\begin{equation}\label{eq:energie}
\begin{aligned}
E^\eps(t):=&\frac{1}{2}\|\eps\nabla_x u^\eps(t)\|^2_{L^2}  + \int_{\R^n}
V(x)|u^\eps 
(t,x)|^2 dx+ \frac{\eps^{n\sigma}}{\sigma
+1}\|u^\eps(t)\|^{2\sigma +2}_{L^{2\sigma +2}}\\
 = & E^\eps(0)=O(1).
\end{aligned}
\end{equation}
The conservation of energy provides an a priori estimate for $H^\eps
u^\eps$ and $J^\eps u^\eps$ thanks to the identity,
\begin{equation}\label{eq:algebre}
\forall t,x,\ \ |H^\eps(t)u^\eps(t,x)|^2 +\eps^2
|J^\eps(t)u^\eps(t,x)|^2 = |x|^2|u^\eps(t,x)|^2 + |\eps\nabla_x u^\eps(t,x)|^2.
\end{equation}
\begin{prop}\label{prop:global}
Fix $\eps \in ]0,1]$ and let $f,r^\eps \in \Sigma$. Then (\ref{eq:nl}) has
a unique solution $u^\eps \in 
Y^\eps(\R)$ and
there exists $C$ such that for any $t\geq 0$ and
any $\eps\in ]0,1]$, 
\begin{equation}\label{eq:bornesnl}
\|u^\eps(t)\|_{L^2}+\|\eps \nabla_x u^\eps (t)\|_{L^2_x} +\|x u^\eps
(t,x)\|_{L^2_x} \leq C. 
\end{equation}
\end{prop}

\subsection{Propagation before the focus}
\label{sec:avt}
Before the focus, we take as an approximate solution the solution of
the linear problem, that is, $v^\eps$ defined by (\ref{eq:lin}).

Notice that from Prop.~\ref{prop:bkw}, we know the asymptotic
behavior of $v^\eps$ before the focus. We prove that in the very same
region, $v^\eps$ is a good approximation of the nonlinear problem.

\begin{prop}\label{prop:avt}
Assume $f,r^\eps\in \Sigma$. Then 
\begin{equation*}
\limsup_{\eps \rightarrow 0} \sup_{0\leq t\leq\frac{\pi}{2}-\Lambda
\eps}\left\|A^\eps(t)\left( u^\eps(t) -
v^\eps 
(t)\right)\right\|_{L^2} \Tend \Lambda {+\infty} 0,
\end{equation*}
where $A^\eps(t)$ is either of the operators $\Id$, $J^\eps(t)$ or
$H^\eps(t)$. 
\end{prop}
\begin{proof}
Define the remainder $w^\eps = u^\eps - v^\eps$. It solves
\begin{equation*}
\left\{
\begin{aligned}
i\eps\d_t w^\eps +\frac{1}{2}\eps^2\Delta w^\eps &
=V(x)w^\eps +\eps^{n\sigma} |u^\eps|^{2\sigma} u^\eps,\\ 
w^\eps_{\mid t=0}& = r^\eps.
\end{aligned}\right.
\end{equation*}
From Duhamel's principle, this writes,
\begin{equation}\label{eq:DuhW}
w^\eps(t)= U^\eps(t)r^\eps -i \eps^{n\sigma -1}\int_0^t
U^\eps(t-s)\left( |u^\eps|^{2\sigma}u^\eps\right)(s)ds. 
\end{equation}
Since $v^\eps$ solves the linear equation (\ref{eq:lin}), so does
$J^\eps(t)v^\eps$, and 
$$\left\|v^\eps(t)\right\|_{L^2}=\|f\|_{L^2},\ \
\|J^\eps(t)v^\eps\|_{L^2}=\|\nabla f\|_{L^2}.$$
From Sobolev inequality (\ref{eq:sobolevn}), 
$$\left\|v^\eps(t)\right\|_{L^{\underline s}}\leq \frac{C}{|\cos t|^{\delta
({\underline s})}} \|f\|_{L^2}^{1-\delta
({\underline s})}\|\nabla f\|_{L^2}^{\delta
({\underline s})}.$$
Therefore there exists $C_0$ such that 
\begin{equation}\label{eq:estv}
\left\|v^\eps(t)\right\|_{L^{\underline s}}\leq \frac{C_0}{|\cos t|^{\delta
({\underline s})}}.
\end{equation}
From Sobolev inequality, for $\eps$
sufficiently small, 
$\left\|w^\eps(0)\right\|_{L^{\underline s}}< C_0$. From
Prop.~\ref{prop:global}, for fixed $\eps >0$, $u^\eps \in
C(\R,\Sigma)$, and the same obviously holds for $v^\eps$. Therefore,
there exists $t^\eps>0$ such that
\begin{equation}\label{eq:solong}
\left\|w^\eps(t)\right\|_{L^{\underline s}}\leq \frac{C_0}{|\cos t|^{\delta
({\underline s})}},
\end{equation}
for any $t\in [0,t^\eps]$. So long as (\ref{eq:solong}) holds, 
we have 
$$\left\|u^\eps(t)\right\|_{L^{\underline s}}\leq \frac{2C_0}{|\cos t|^{\delta
({\underline s})}},$$
and we can apply Prop.~\ref{prop:estgen}. Indeed, take ${\EuScript
U}^\eps = U^\eps$, 
$h^\eps=\eps^{n\sigma}|u^\eps|^{2\sigma}v^\eps$ 
and $F^\eps(w^\eps)= |u^\eps|^{2\sigma}w^\eps$. From H\"older inequality
and the above estimate,
\begin{equation*}
\begin{aligned}
\left\|F^\eps(w^\eps)(t)\right\|_{L^{{\underline r}'}}& \leq
\|u^\eps(t)\|^{2\sigma}_{L^{\underline{s}}}
\|w^\eps(t)\|_{L^{\underline r}}\\
&\leq \frac{(2C_0)^{2\sigma}}{\left(|\cos 
t| \right)^{2\sigma
\delta(\underline{s})}}\|w^\eps(t)\|_{L^{\underline r}}.
\end{aligned}
\end{equation*}
Assume (\ref{eq:solong}) holds for $0\leq t\leq T$. If $0 \leq t\leq
T\leq \frac{\pi}{2}-\Lambda \eps$, then $\eps \lesssim \cos 
t$, and the 
above estimate shows that $F^\eps$ satisfies assumption (\ref{eq:hypF}). 

From Cor.~\ref{cor:estgen}, if $\Lambda$ is sufficiently large, then
for $0\leq t\leq T\leq \frac{\pi}{2}-\Lambda \eps$, and for any
$(q,r)$ admissible, 
\begin{equation*}\label{eq:wL2}
\|w^\eps\|_{L^\infty(0,T;L^2)}\leq C\|r^\eps\|_{L^2}+
C\eps^{n\sigma-1-\frac{1}{q}} \left\||u^\eps|^{2\sigma}v^\eps\right\|_
{ L^{q'}(0,T;L^{r'})}.  
\end{equation*}
Taking $(q,r)=(\underline q, \underline r)$ yields, from H\"older
inequality,
$$\left\||u^\eps|^{2\sigma}v^\eps\right\|_
{ L^{{\underline q}'}(0,T;L^{{\underline r}'})}\leq
\|u^\eps\|^{2\sigma}_{L^{\underline k}(0,T;L^{\underline s})}\|v^\eps\|_
{ L^{\underline q}(0,T;L^{\underline r})}.$$
The first term of the right-hand side is estimated through (\ref{eq:estv})
and (\ref{eq:solong}). The last term is estimated the same way, for
(\ref{eq:estv}) still holds when replacing $\underline s$ with
$\underline r$. Therefore,
$$\left\||u^\eps|^{2\sigma}v^\eps\right\|_
{ L^{{\underline q}'}(0,T;L^{{\underline r}'})}\leq \frac{C}{\left(
\frac{\pi}{2}-T \right)^{n\sigma -1 -\frac{1}{\underline q}}},$$
and 
\begin{equation}\label{eq:reste}
\|w^\eps\|_{L^\infty(0,T;L^2)}\leq C\|r^\eps\|_{L^2}+ C 
\left(\frac{\eps}{\frac{\pi}{2}-T} \right)^{n\sigma -1
-\frac{1}{\underline q}}. 
\end{equation}

Now apply the operator $J^\eps$ to (\ref{eq:DuhW}). Since $J^\eps$ and
$U^\eps$ commute, it yields,
\begin{equation*}
J^\eps(t)w^\eps= U^\eps(t)J^\eps(0)r^\eps -i \eps^{n\sigma -1}\int_0^t
U^\eps(t-s)J^\eps(s)\left( |u^\eps|^{2\sigma}u^\eps\right)(s)ds. 
\end{equation*}
Because $J^\eps$ acts on this nonlinear like a derivative, we have an
equation which is very similar to (\ref{eq:DuhW}), with $w^\eps$
replaced by $J^\eps w^\eps$ and $r^\eps$ replaced by $-i\nabla
r^\eps$. Therefore the same computation as above yields
\begin{equation}\label{eq:reste2}
\|J^\eps w^\eps\|_{L^\infty(0,T;L^2)}\leq C\|\nabla r^\eps\|_{L^2}+ C 
\left(\frac{\eps}{\frac{\pi}{2}-T} \right)^{n\sigma -1
-\frac{1}{\underline q}}. 
\end{equation} 
Combining (\ref{eq:reste}) and (\ref{eq:reste2}) yields, along with
(\ref{eq:sobolevn}),
$$ \forall t\in [0,T],\ \  \left\|w^\eps(t)\right\|_{L^{\underline
s}}\leq \frac{C}{|\cos t|^{\delta 
({\underline s})}}\left(\|r^\eps\|_{H^1} +
\left(\frac{\eps}{\frac{\pi}{2}-t} \right)^{n\sigma -1 
-\frac{1}{\underline q}} \right).$$
Therefore, choosing $\eps$ sufficiently small and $\Lambda$
sufficiently large, we deduce that we can take
$T=\frac{\pi}{2}-\Lambda \eps$. This yields Prop.~\ref{prop:avt} for
$A^\eps=\Id$ and $J^\eps$. The case $A^\eps=H^\eps$ is now
straightforward.   
\end{proof}

\subsection{Matching linear and nonlinear r\'egimes}
\label{sec:matching}
When time approaches $\pi/2$, the nonlinear term cannot be neglected. On
the other hand, since the solution tends to concentrate at the origin,
the potential becomes negligible. It is then natural to seek an
approximate solution $\tilde v^\eps$ that solves
$$i\eps \d_t \tilde v^\eps +\frac{1}{2}\eps^2 \Delta \tilde v^\eps =
\eps^{n\sigma} |\tilde v^\eps|^{2\sigma}\tilde v^\eps.$$
The question that arises naturally is, how can we match $\tilde
v^\eps$ and $v^\eps$? With the results of \cite{Ca2} in mind, we can
expect that $\tilde v^\eps$ is exactly a concentrating profile,
\begin{equation}\label{eq:scaling}
\tilde v^\eps(t,x)= \frac{1}{ \eps^{n/2}}\psi \left( 
\frac{t-\frac{\pi}{2}}{\eps}, \frac{x}{\eps}\right).
\end{equation}
The function $\psi$ must be defined to match the solution $u^\eps$, or
one of 
its approximations $v^\eps$ or $v^\eps_{\rm app}$, when
$t=\pi/2-\Lambda\eps$, for $\Lambda$ 
sufficiently large. Notice that this problem was already encountered
by H.~Bahouri and P.~G\'erard in \cite{BG3} (see also \cite{BG},
I.~Gallagher and P.~G\'erard \cite{GG}). We
prove that for $\Lambda >0$ sufficiently large, the propagation for
$\pi/2-\Lambda \eps\leq t\leq \pi/2+\Lambda \eps$ is described by
$\tilde v^\eps$. 

Write $t_*^\eps =\pi/2-\Lambda\eps$, and assume from now on that
$\Lambda >1$. For large $\Lambda$,
Prop.~\ref{prop:bkw} and \ref{prop:avt} imply 
$$u^\eps(t_*^\eps,x)\sim v^\eps(t_*^\eps,x) \sim v^\eps_{\rm
app}(t_*^\eps,x)=\frac{1}{({\sin (\Lambda \eps)})^{n/2}}
f\left(\frac{x}{\sin (\Lambda \eps)} \right)e^{-i\frac{|x|^2}{2\eps\tan
(\Lambda\eps) }}.$$
For $\Lambda \eps$ close to zero, the following approximation is
expected,
$$\frac{1}{(\sin (\Lambda \eps))^{n/2}}
f\left(\frac{x}{\sin (\Lambda \eps)} \right)e^{-i\frac{|x|^2}{2\eps\tan
(\Lambda\eps) }}\sim \frac{1}{(\Lambda \eps)^{n/2}}
f\left(\frac{x}{\Lambda \eps} \right)e^{-i\frac{|x|^2}{2\eps
(\Lambda\eps) }}.$$
We prove that this approximation is correct in
Lemma~\ref{lem:1} below. From (\ref{eq:scaling}), this should also be close
to
$$\frac{1}{ \eps^{n/2}}\psi\left(-\Lambda,\frac{x}{\eps} \right).$$
Recall the classical result,
\begin{prop}\label{prop:scattnls}(\cite{HT87}, Theorem~1.1; \cite{CW92},
Theorem~4.2) Assume $\psi_- \in \Sigma$ and
$\frac{2}{n+2}<\sigma<\frac{2}{n-2}$ if $n\geq 2$, $\sigma >1$ if
$n=1$. Denote 
$$\sigma_0(n):=\frac{2-n +\sqrt{n^2+12n +4}}{4n}.$$
If
$\sigma>\sigma_0(n)$ or if $\|\psi_-\|_\Sigma$ is
sufficiently small, then 
\begin{itemize}
\item There exists a unique $\psi \in C(\R_t,\Sigma)$ solution of
(\ref{eq:nls}), such that 
$$\lim_{t\rightarrow -\infty}\|\psi_- -U_0(-t)\psi(t)\|_\Sigma
=0,\textrm{ \ \ where } U_0(t)=e^{i\frac{t}{2}\Delta}.$$
\item There exists a unique $\psi_+ \in \Sigma$ such that 
$$\lim_{t\rightarrow +\infty}\|\psi_+ -U_0(-t)\psi(t)\|_\Sigma =0.$$
\end{itemize}
\end{prop}
Recall that the asymptotic state $\psi_-$ was defined in introduction
by, 
\begin{equation*}
\psi_-:=\frac{1}{(2i\pi)^{n/2}}\widehat f,
\end{equation*}
and the approximate solution (near $t=\pi/2$) is given by 
\begin{equation*}
\tilde v^\eps(t,x)= \frac{1}{\eps^{n/2}}\psi \left( 
\frac{t-\frac{\pi}{2}}{\eps}, \frac{x}{\eps}\right).
\end{equation*}
We prove, 
\begin{prop}\label{prop:matching}
Assume $f,r^\eps\in \Sigma$. 
Take $\psi_-$ defined by (\ref{eq:psi-}). Then 
\begin{equation*}
\limsup_{\eps \rightarrow 0}
\left\| u^\eps\left(\frac{\pi}{2}-\Lambda \eps, .\right)
-\frac{1}{\eps^{n/2}}\left(U_0(-\Lambda)\psi_-\right)
\left(\frac{.}{\eps}\right)\right\|_{L^2}\Tend \Lambda {+\infty} 0,
\end{equation*}
and the same holds when applying either of the operators
$\eps\nabla_x$ or $\frac{x}{\eps} -i\Lambda\eps\nabla_x$ to the considered
functions. 
\end{prop}
\begin{proof}
From Cor.~\ref{cor:idem} (from which $v^\eps \sim
v^\eps_{\rm app}$)
and Prop.~\ref{prop:avt} (from which $v^\eps \sim u^\eps$),  
\begin{equation}\label{eq:step1}
\limsup_{\eps \rightarrow 0}
\left\| u^\eps(t_*^\eps,x)-\frac{1}{(\sin (\Lambda \eps))^{n/2}}
f\left(\frac{x}{\sin (\Lambda \eps)} \right)e^{-i\frac{|x|^2}{2\eps\tan
(\Lambda\eps) }}\right\|_{L^2}\Tend \Lambda {+\infty} 0, 
\end{equation}
and the same result holds when applying either of the operators
$J^\eps(t_*^\eps)$ or $H^\eps(t_*^\eps)$. Notice that applying
$J^\eps(t_*^\eps)$ or $H^\eps(t_*^\eps)$ 
is not so different from applying $\eps\nabla_x$ or $\frac{x}{\eps}
-i\Lambda\eps\nabla_x$, for when $\Lambda \eps$ goes to zero,
$$J^\eps(t_*^\eps) \sim \frac{x}{\eps} -i\Lambda\eps\nabla_x\ \ \ \ ,\ 
H^\eps(t_*^\eps)\sim i\eps\nabla_x.$$
Recall that $t_*^\eps=\pi/2 -\Lambda \eps$.   

\begin{lem}\label{lem:2}
Let $a^\eps(t_*^\eps,.)\in \Sigma$ be a family of functions 
such that there exists $C_*$ independent of $\eps\in ]0,1]$ such that,
\begin{equation}\label{eq:assa}
\|xa^\eps(t_*^\eps,x)\|_{L^2} +\left\| \eps \nabla_x
a^\eps(t_*^\eps,x)\right\|_{L^2} \leq C_*.
\end{equation}
Then for any $\Lambda >1$, 
\begin{equation*}
\begin{aligned}
\limsup_{\eps \rightarrow 0}\left\|\left(J^\eps(t_*^\eps)
-\frac{x}{\eps}+i\Lambda \eps \nabla_x \right)
a^\eps(t_*^\eps)\right\|_{L^2}&=\\
=\limsup_{\eps \rightarrow 0}\|(H^\eps(t_*^\eps)
&-i \eps \nabla_x ) a^\eps(t_*^\eps)\|_{L^2}=0.
\end{aligned}
\end{equation*}
In particular, we can take  $a^\eps=u^\eps$ or $a^\eps=v^\eps_{\rm app}$. 
\end{lem}
\rq Lemma \ref{lem:2} has a simple geometric interpretation. Near the
focus, rays of geometric optics, given by (\ref{eq:rays}), are
straightened (Fig.~\ref{fig:zoom}). 
\begin{figure}[htbp]
\begin{center}
\input{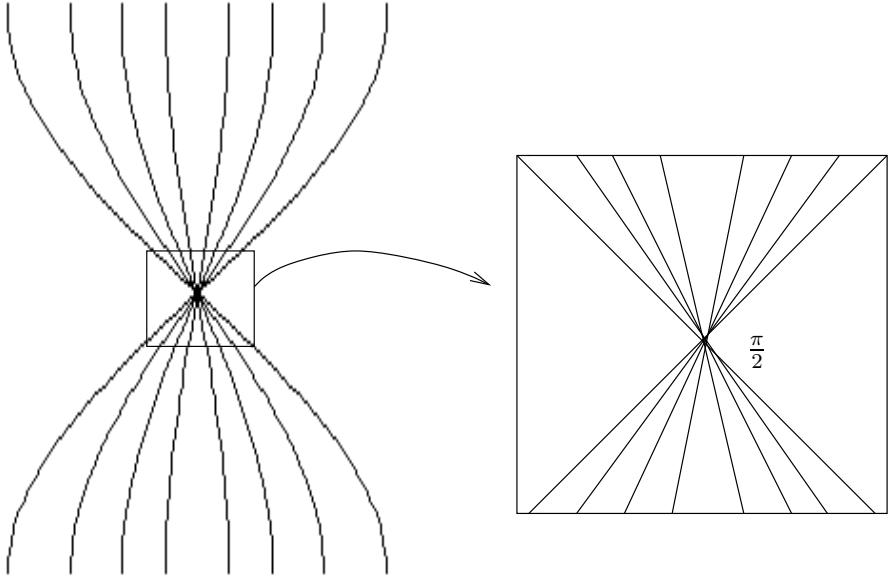}
\caption{Rays of geometric optics are straightened near $t=\frac{\pi}{2}$.}
\label{fig:zoom}
\end{center}
\end{figure} 
Thus in the neighborhood of
$t=\pi/2$, rays are almost 
straight lines, that is, the geometry is nearly the same as in
\cite{Ca2}. In that case, with the natural scaling (\ref{eq:scaling}),
the ``good'' operators are $\eps\nabla_x$ and
$\frac{x}{\eps}+i(t-\pi/2)\nabla_x$. \\

\begin{proof}[Proof of Lemma~\ref{lem:2}]
 Fix $\Lambda >1$. 
\begin{equation*}
\left(J^\eps(t_*^\eps)
-\frac{x}{\eps}+i\Lambda \eps \nabla_x \right)
a^\eps(t_*^\eps,x)= \left( (\cos (\Lambda \eps)-1)\frac{x}{\eps} -i(\sin
(\Lambda \eps) -\Lambda \eps)\nabla_x\right)a^\eps(t_*^\eps,x).
\end{equation*}
Taking the $L^2$ norm yields,
$$\left\|\left(J^\eps(t_*^\eps)
-\frac{x}{\eps}+i\Lambda \eps \nabla_x \right)
a^\eps(t_*^\eps)\right\|_{L^2} \leq C(\Lambda \eps)^2 \left\|
\frac{x}{\eps} a^\eps(t_*^\eps,x)\right\|_{L^2_x} + C(\Lambda \eps)^3
\|\nabla_x a^\eps(t_*^\eps)\|_{L^2_x}.$$
The assumption (\ref{eq:assa}) (which is a consequence of (\ref{eq:bornesnl})
for $u^\eps$, and straightforward for $v^\eps_{\rm app}$) implies
$$\left\|\left(J^\eps(t_*^\eps)
-\frac{x}{\eps}+i\Lambda \eps \nabla_x \right)
a^\eps(t_*^\eps)\right\|_{L^2} \leq C\Lambda^2\eps + C
\Lambda^3\eps^2, $$
which proves the first part of the lemma. Similarly,
\begin{equation*}
\begin{aligned}
\left\|\left(H^\eps(t_*^\eps)
-i \eps \nabla_x \right)
a^\eps(t_*^\eps)\right\|_{L^2} &\leq C(\Lambda \eps) \left\|
x a^\eps(t_*^\eps,x)\right\|_{L^2_x} + C(\Lambda \eps)^2
\|\eps \nabla_x a^\eps(t_*^\eps)\|_{L^2_x}\\
&\leq C(\Lambda \eps) + C(\Lambda \eps)^2.
\end{aligned}
\end{equation*}
This completes the proof of the lemma.  
\end{proof}

Now we prove that in \eqref{eq:step1}, we can replace $\sin (\Lambda
\eps)$ and $\tan(\Lambda\eps)$ with $\Lambda\eps$ up to a small
error term. Denote 
$$ \tilde v^\eps_{\rm app}(t,x)=\frac{1}{\left(\frac{\pi}{2}-
t\right)^{n/2}}f\left(\frac{x}{\frac{\pi}{2}- t} 
\right) e^{-i\frac{|x|^2}{2\eps(\pi/2-t)}}.$$
\begin{lem}\label{lem:1}
Assume $f\in \Sigma$. For any $\Lambda >1$,
\begin{equation*}
\begin{aligned}
\limsup_{\eps \rightarrow 0}\left\| (v^\eps_{\rm app}-\tilde
v^\eps_{\rm app})(t_*^\eps) 
\right\|_{L^2} =  
\limsup_{\eps \rightarrow 0}\left\| \eps\nabla_x (v^\eps_{\rm app}-\tilde
v^\eps_{\rm app})(t_*^\eps) \right\|_{L^2}=\\ 
=\limsup_{\eps \rightarrow 0}\left\| 
\left( \frac{x}{\eps} -i\Lambda\eps\nabla_x\right)(v^\eps_{\rm app}-
\tilde v^\eps_{\rm app})(t_*^\eps) \right\|_{L^2}
&=0.
\end{aligned}
\end{equation*}
\end{lem}
\begin{proof}[Proof of Lemma~\ref{lem:1}]
  Write $\lambda =\Lambda
\eps$. For fixed $\Lambda$, $\lambda$ is a small parameter when $\eps$
goes to zero, and
\begin{equation*}
\begin{aligned}
(v^\eps_{\rm app}-\tilde v^\eps_{\rm app})(t_*^\eps,x)=& 
\frac{1}{(\sin \lambda)^{n/2}}
f\left(\frac{x}{\sin \lambda} \right)e^{-i\frac{|x|^2}{2\eps\tan
\lambda }} -\frac{1}{\lambda^{n/2}}
f\left(\frac{x}{\lambda} \right)e^{-i\frac{|x|^2}{2\eps
\lambda }} \\
=& \left( \frac{1}{(\sin \lambda)^{n/2}}
f\left(\frac{x}{\sin \lambda} \right)-\frac{1}{\lambda^{n/2}}
f\left(\frac{x}{ \lambda} \right)\right)e^{-i\frac{|x|^2}{2\eps\tan
\lambda }}\\
& +\frac{1}{\lambda^{n/2}}
f\left(\frac{x}{ \lambda} \right)\left( e^{-i\frac{|x|^2}{2\eps\tan
\lambda }}-e^{-i\frac{|x|^2}{2\eps
\lambda }}\right).
\end{aligned}
\end{equation*}
Taking the $L^2$ norm yields,
\begin{equation*}
\begin{aligned}
\|(v^\eps_{\rm app}-\tilde v^\eps_{\rm app})(t_*^\eps)\|_{L^2} \leq &
\left\| \frac{1}{(\sin \lambda)^{n/2}}
f\left(\frac{.}{\sin \lambda} \right)-\frac{1}{\lambda^{n/2}}
f\left(\frac{.}{ \lambda} \right)\right\|_{L^2}\\
& +  \left\|f(x)\left( e^{-i\frac{\lambda^2|x|^2}{2\eps\tan
\lambda }}-e^{-i\frac{\lambda |x|^2}{2\eps
}}\right)\ \right\|_{L^2}\\
\leq & \left\| \left(\frac{1}{(\sin \lambda)^{n/2}}
-\frac{1}{\lambda^{n/2}}\right)
f\left(\frac{.}{ \sin \lambda} \right)\right\|_{L^2} \\
&+
\left\| \frac{1}{\lambda^{n/2}}
\left(f\left(\frac{.}{\sin \lambda} \right)-
f\left(\frac{.}{ \lambda} \right)\right)\right\|_{L^2}\\
& +  \left\|f(x)\left( e^{-i\frac{|x|^2}{2\eps}\left( \frac{\lambda^2}{\tan
\lambda}  -\lambda\right)}-1\right)\ \right\|_{L^2}\\
\leq & \left| \left(\frac{\sin \lambda}{\lambda}\right)^{n/2}-1\right|
\|f\|_{L^2}+ \left\| 
f\left(\frac{\lambda \ .}{\sin \lambda} \right)-
f(.)\right\|_{L^2}  \\
&+ \left\|f(x)\left( e^{-i\frac{|x|^2}{2\eps}\left( \frac{\lambda^2}{\tan
\lambda}  -\lambda\right)}-1\right)\ \right\|_{L^2}.
\end{aligned}
\end{equation*}
The first term of the right-hand side clearly goes to zero with
$\lambda$. So does the second one: if $f\in C_0^\infty(\R)$, it is
$O(\lambda^2)$, and by density, it is $o(1)$ when $\lambda$ goes to
zero for any $f\in L^2$. 
Recalling that $\lambda =\Lambda \eps$, we have
$$\frac{1}{\eps}\left( \frac{\lambda^2}{\tan 
\lambda}  -\lambda\right) = \Lambda \left( \frac{\lambda}{\tan 
\lambda}  -1\right).$$
Thus, for any {\it fixed} $\Lambda >1$, this term goes to zero when
$\eps$ goes to zero. Therefore, from dominated convergence,  for any
fixed $\Lambda >1$,
$$\limsup_{\eps \rightarrow 0}\left\|f(x)\left
( e^{-i\frac{|x|^2}{2\eps}\left( \frac{\lambda^2}{\tan 
\lambda}  -\lambda\right)}-1\right)\ \right\|_{L^2} =0.$$
Computations for $\left\| \eps\nabla_x (v^\eps_{\rm app}-\tilde
v^\eps_{\rm app})(t_*^\eps) \right\|_{L^2}$ and 
$\left\| 
\left( \frac{x}{\eps} -i\Lambda\eps\nabla_x\right)(v^\eps_{\rm app}-
\tilde v^\eps_{\rm app})(t_*^\eps) \right\|_{L^2}$ are similar  and
essentially involve one more derivative or one more momentum. Indeed,
$v^\eps_{\rm app}$ and $\tilde v^\eps_{\rm app}$ behave well with
respect to the operators $\eps\nabla_x$  and
$\frac{x}{\eps} -i\Lambda\eps\nabla_x$, thus we can use the same
density argument as above. 
\end{proof}

The next step to prove Prop.~\ref{prop:matching} consists in
comparing $\tilde v^\eps_{\rm app}$ and the rescaled free evolution of
the asymptotic state $\psi_-$. 

\begin{lem}\label{lem:3}
Assume $f\in \Sigma$. 
The following limits hold, uniformly with respect to $\eps \in ]0,1]$,
\begin{equation*}
\begin{aligned}
\lim_{\Lambda \rightarrow +\infty}
& \left\|\tilde v^\eps_{\rm app} (t_*^\eps)
- \frac{1}{ \eps^{n/2}} \left(U_0(-\Lambda)\psi_- 
\right)\left( \frac{.}{\eps}\right)\right\|_{L^2}=\\
=& \lim_{\Lambda \rightarrow +\infty}
 \left\|\eps\nabla_x \left(\tilde v^\eps_{\rm app}
(t_*^\eps) 
- \frac{1}{\eps^{n/2}} \left(U_0(-\Lambda)\psi_- 
\right)\left( \frac{.}{\eps}\right)\right)\right\|_{L^2}\\
=&\lim_{\Lambda \rightarrow +\infty}
 \left\|\left( \frac{x}{\eps}-i\Lambda
 \eps \nabla_x \right)
 \left(\tilde v^\eps_{\rm app}
(t_*^\eps) 
- \frac{1}{\eps^{n/2}} \left(U_0(-\Lambda)\psi_- 
\right)\left( \frac{.}{\eps}\right)\right)\right\|_{L^2}=0.
\end{aligned}
\end{equation*}
\end{lem}
\begin{proof}[Proof of Lemma~\ref{lem:3}]
 From the Fourier Inversion
 Formula, we have
\begin{equation*}
  \tilde v^\eps_{\rm app} (t_*^\eps,x)=\frac{1}{(\Lambda
  \eps)^{n/2}}f\left(\frac{x}{\Lambda \eps}  
\right) e^{-i\frac{|x|^2}{2\eps^2\Lambda}}
=\frac{1}{(2\pi)^n}
\frac{1}{(\Lambda \eps)^{n/2}} e^{-i\frac{|x|^2}{2\eps^2\Lambda}} \int
e^{i \frac{x.y}{\eps\Lambda}}\widehat f (y)dy . 
\end{equation*}
On the other hand, the expression of the free Schr\"odinger group
$U_0$ implies, along with definition (\ref{eq:psi-}), 
\begin{equation*}
\begin{aligned}
\frac{1}{\eps^{n/2}} \left(U_0(-\Lambda)\psi_- 
\right)\left( \frac{x}{\eps}\right)&=
\left(\frac{i}{2\pi\Lambda \eps}\right)^{n/2}
 e^{-i\frac{|x|^2}{2\eps^2\Lambda}} \int
e^{i \frac{x\cdot y}{\eps\Lambda}-i\frac{|y|^2}{2\Lambda}}\psi_- (y)dy\\
&=\frac{1}{(2\pi)^n}
\frac{1}{(\Lambda \eps)^{n/2}} e^{-i\frac{|x|^2}{2\eps^2\Lambda}} \int
e^{i \frac{x\cdot y}{\eps\Lambda}-i\frac{|y|^2}{2\Lambda}}\widehat f (y)dy.
\end{aligned}
\end{equation*}
Thus the remainder we have to assess writes
$$\frac{1}{(2\pi)^n}
\frac{1}{(\Lambda \eps)^{n/2}} e^{-i\frac{|x|^2}{2\eps^2\Lambda}} \int
e^{i \frac{x\cdot y}{\eps\Lambda}}\left(1 -
e^{-i\frac{|y|^2}{2\Lambda}}\right)\widehat f (y)dy,$$
which is also,
$$ \frac{1}{(\Lambda \eps)^{n/2}} e^{-i\frac{|x|^2}{2\eps^2\Lambda}} 
\left( \left( 1- e^{i\frac{\Delta}{2\Lambda}}\right)f\right)\left
( \frac{x}{\Lambda\eps}\right). $$ 
The lemma then follows from the strong convergence in $L^2$,
$e^{i\delta \Delta} \Tend \delta 0 1$.  
\end{proof}

Lemmas~\ref{lem:2}, \ref{lem:1} and \ref{lem:3} imply
Prop.~\ref{prop:matching}. 
\end{proof}

\subsection{Description of the solution near the focus}
\label{sec:foc} 
Propositions \ref{prop:scattnls} and \ref{prop:matching} imply that
\begin{equation*}
\begin{aligned}
 \limsup_{\eps \rightarrow 0}&
\left\| u^\eps\left(\frac{\pi}{2}-\Lambda \eps, \cdot \right)
-\tilde v^\eps\left(\frac{\pi}{2}-\Lambda \eps,
 \cdot \right)\right\|_{L^2}+\\
&+\limsup_{\eps \rightarrow 0}
\left\| \eps \nabla_x\left(u^\eps\left(\frac{\pi}{2}-\Lambda \eps, \cdot\right)
-\tilde v^\eps\left(\frac{\pi}{2}-\Lambda \eps,
 \cdot \right)\right)\right\|_{L^2}+\\
&+\limsup_{\eps \rightarrow 0}
\left\| \left(\frac{x}{\eps}-i\Lambda\eps\nabla_x
 \right)\left(u^\eps\left(\frac{\pi}{2}-\Lambda \eps, \cdot\right) 
-\tilde v^\eps\left(\frac{\pi}{2}-\Lambda \eps,
 \cdot\right)\right)\right\|_{L^2}
\Tend \Lambda {+\infty} 0.
\end{aligned}
\end{equation*}
This means that taking $\Lambda$ large enough, and $\eps$ small 
enough, the difference $u^\eps- \tilde v^\eps$ is small at time
$t_*^\eps = \pi/2 -\Lambda\eps$, which is the ``initial'' time in the
boundary layer where nonlinear effects take place (and where the
potential is negligible). Since the role of $r^\eps$ is negligible, we
first assume $r^\eps\equiv 0$.

\begin{prop}\label{prop:traversee1}
Assume $f\in {\EuScript H}$, and that the nonlinearity is $C^2$, that
is, $\sigma>1/2$, which is possible only if $n\leq 5$. Then 
the difference $u^\eps- \tilde v^\eps$ is small around the focus.
\begin{equation*}
\limsup_{\eps \rightarrow 0} \sup_{\frac{\pi}{2}-\Lambda \eps \leq
t\leq \frac{\pi}{2}+\Lambda \eps}\left\|A^\eps(t)\left( u^\eps(t) -
\tilde v^\eps 
(t)\right)\right\|_{L^2} \Tend \Lambda {+\infty} 0,
\end{equation*}
where $A^\eps(t)$ is either of the operators $\Id$, $J^\eps(t)$ or
$H^\eps(t)$. 
\end{prop}
\rq The assumption $\sigma >\frac{1}{2}$ is needed to prove
Lemma~\ref{lem:mulder} below. It seems purely technical, and one
expects Lemma~\ref{lem:mulder} to hold without this assumption. If
$n=2$, the nonlinearity is automatically $C^2$ thanks to the
assumption $\sigma >\frac{2}{n+2}$. If $n=3$, then we have to restrict
our study to the case $\frac{1}{2}<\sigma< 2$. In particular, the
value $\sigma =1$, which corresponds to a cubic nonlinearity, is
accepted. 

\begin{proof}
Propositions \ref{prop:scattnls} and \ref{prop:matching} imply that
\begin{equation*}
\begin{aligned}
 \limsup_{\eps \rightarrow 0}&
\left\| u^\eps\left(\frac{\pi}{2}-\Lambda \eps, \cdot\right)
-\tilde v^\eps\left(\frac{\pi}{2}-\Lambda \eps,
 \cdot\right)\right\|_{L^2}+\\
&+\limsup_{\eps \rightarrow 0}
\left\| \eps \nabla\left(u^\eps\left(\frac{\pi}{2}-\Lambda \eps, \cdot\right)
-\tilde v^\eps\left(\frac{\pi}{2}-\Lambda \eps,
 \cdot\right)\right)\right\|_{L^2}+\\
&+\limsup_{\eps \rightarrow 0}
\left\| \left(\frac{x}{\eps}-i\Lambda\eps\nabla
 \right)\left(u^\eps\left(\frac{\pi}{2}-\Lambda \eps, \cdot\right) 
-\tilde v^\eps\left(\frac{\pi}{2}-\Lambda \eps,
 \cdot\right)\right)\right\|_{L^2}
\Tend \Lambda {+\infty} 0.
\end{aligned}
\end{equation*}

Define the remainder $\tilde w^\eps =  u^\eps-\tilde v^\eps $,
and keep the notation $t_*^\eps = \pi/2 -\Lambda\eps$. From
Prop.~\ref{prop:matching}, 
\begin{equation*}
 \limsup_{\eps \rightarrow 0}
\left\| B^\eps(t_*^\eps)\tilde w^\eps(t_*^\eps)\right\|_{L^2}
\Tend \Lambda {+\infty} 0,
\end{equation*}
where $B^\eps(t)$ is either of the operators $\Id$,
$\frac{x}{\eps}+i(t-\pi/2) \nabla$ or
$\eps\nabla$. From Lemma~\ref{lem:2}, this implies 
\begin{equation*}
 \limsup_{\eps \rightarrow 0}
\left\| A^\eps(t_*^\eps)\tilde w^\eps(t_*^\eps)\right\|_{L^2}
\Tend \Lambda {+\infty} 0,
\end{equation*}
where $A^\eps(t)$ is either of the operators $\Id$, $J^\eps(t)$ or
$H^\eps(t)$. 

From the conservation of energy (\ref{eq:energie}), we have
$$\|\eps\nabla u^\eps(t)\|_{L^2}\leq C.$$
From the conservation of energy for (\ref{eq:nls}), 
$$\frac{d}{dt}\left( \frac{1}{2}
\|\nabla \psi(t)\|_{L^2}^2 + \frac{1}{\sigma +1} 
\|\psi(t)\|_{L^{2\sigma +2}}^{2\sigma +2} \right)=0,$$
we have 
$$\|\eps\nabla\tilde v^\eps(t)\|_{L^2}\leq C.$$
Therefore, since $\tilde w^\eps =u^\eps- \tilde v^\eps$,
$$\|\eps\nabla \tilde w^\eps(t)\|_{L^2}\leq C.$$
From Sobolev inequality,
$$\|\tilde w^\eps(t)\|_{L^{\underline s}} \leq C \|\tilde
w^\eps(t)\|_{L^2}^{1-\delta(\underline s)} \|\nabla\tilde
 w^\eps(t)\|_{L^2}^{\delta(\underline s)},$$ 
and there exists $C_1$ such that for any $t \in \R$, 
\begin{equation}\label{eq:solong2}
\|\tilde w^\eps(t)\|_{L^{\underline s}} \leq \frac{C_1}{
\eps^{\delta(\underline s)}}.
\end{equation}
This estimate will be useful for $|t-\pi/2|\leq \Lambda_0 \eps$, where
$\Lambda_0$ is given by Cor.~\ref{cor:estgen}. For $|t-\pi/2|\geq
\Lambda_0 \eps$, sharper estimates are provided by $J^\eps$, along with
Sobolev  inequality (\ref{eq:sobolevn}).

The first step of the proof consists in showing that the harmonic
potential can be truncated near the origin without altering the
asymptotics. Let $\chi\in C_0^\infty (\R^n)$ be a cut-off function, with
$$\operatorname{supp}\chi \subset B(0,2), \ \ \ 0\leq \chi \leq 1
\textrm{ and } \  \forall x\in B(0,1),\ 
\chi(x)=1.$$
For $R>0$, define 
$$u_R^\eps(t,x)= \chi \left( \frac{x}{R}\right)u^\eps(t,x).$$
\begin{lem}\label{lem:troncature}
Assume $f\in {\EuScript H}$, $\sigma > \frac{1}{2}$  and 
take $R=\eps^\alpha$. Then for any $0<\alpha <1$, 
\begin{equation*}
\limsup_{\eps \rightarrow 0} \sup_{\frac{\pi}{2} -\Lambda\eps\leq t 
\leq \frac{\pi}{2} +\Lambda\eps}  \|A^\eps(t)(u^\eps(t)
- u_R^\eps(t))\|_{L^2}\Tend \Lambda {+ \infty} 0,
\end{equation*}
where $A^\eps(t)$ is either of the operators $\Id$, $J^\eps(t)$ or
$H^\eps(t)$.
\end{lem}
\begin{proof}[Proof of Lemma~\ref{lem:troncature}]
The function
$u_R^\eps$ satisfies,
\begin{equation*}
\left(i\eps \d_t +\frac{1}{2}\eps^2\Delta -\frac{|x|^2}{2} \right)u_R^\eps=
\eps^{n\sigma} |u^\eps|^{2\sigma}u_R^\eps+\frac{\eps^2}{2R}\nabla\chi\left(
\frac{x}{R}\right)\cdot\nabla u^\eps +\left( \frac{\eps}{R}\right)^2
\Delta\chi\left( \frac{x}{R}\right) u^\eps,
\end{equation*}
therefore the difference $w_R^\eps :=u^\eps- u_R^\eps$ solves,
\begin{equation*}
\left(i\eps \d_t +\frac{1}{2}\eps^2\Delta -\frac{|x|^2}{2} \right)w_R^\eps=
\eps^{n\sigma} |u^\eps|^{2\sigma}w_R^\eps-\frac{\eps^2}{2R}\nabla\chi\left(
\frac{x}{R}\right)\nabla u^\eps -\left( \frac{\eps}{R}\right)^2
\Delta\chi\left( \frac{x}{R}\right) u^\eps.
\end{equation*}
From Lemma~\ref{lem:estimL2}, and because the term $\eps^{n\sigma}
|u^\eps|^{2\sigma}$ can be considered as a real potential, 
$$\eps\d_t \| w_R^\eps(t)\|_{L^2}\leq C \frac{\eps}{R}\|\eps\nabla
u^\eps(t)\|_{L^2}  + C\left( \frac{\eps}{R}\right)^2\|
u^\eps(t)\|_{L^2},$$
which implies, from (\ref{eq:masse}) and (\ref{eq:bornesnl}), 
$$\eps\d_t \| w_R^\eps(t)\|_{L^2}\leq C \frac{\eps}{R}
  + C\left( \frac{\eps}{R}\right)^2.$$
Integrating this inequality on $\left[ \frac{\pi}{2} -\Lambda\eps, 
\frac{\pi}{2} +\Lambda\eps \right]$ gives
$$\sup_{\frac{\pi}{2} -\Lambda\eps\leq t 
\leq \frac{\pi}{2} +\Lambda\eps}  \|w_R^\eps(t)\|_{L^2} \leq 
\|w_R^\eps(\pi/2 -\Lambda \eps)\|_{L^2}
+C \Lambda\frac{\eps}{R}
  + C\Lambda\left( \frac{\eps}{R}\right)^2.$$
Taking $R=\eps^\alpha$ with $0<\alpha <1$ yields,
$$\limsup_{\eps \rightarrow 0}\sup_{\frac{\pi}{2} -\Lambda\eps\leq t 
\leq \frac{\pi}{2} +\Lambda\eps}  \|w_R^\eps(t)\|_{L^2} \leq 
\limsup_{\eps \rightarrow 0}\|w_R^\eps(\pi/2 -\Lambda \eps)\|_{L^2}.$$
Now since $\psi_- \in L^2$, $0<\alpha <1$ implies, along with the
dominated convergence theorem,
$$ \left\| \left(1-\chi
\left(\frac{\cdot}{\eps^\alpha}\right)\right)\frac{1}{\eps^{n/2}}
\left(U_0(-\Lambda)\psi_-\right)  
\left(\frac{.}{\eps}\right)\right\|_{L^2}\Tend \eps 0 0.$$
From Prop.~\ref{prop:matching}, the first part of
Lemma~\ref{lem:troncature} (with $A^\eps =Id$) follows. \\

To estimate $J^\eps w_R^\eps$, notice that 
$$J^\eps(t) w_R^\eps(t,x) = \left(1-\chi \left(\frac{x}{R}\right)
\right)J^\eps(t) u^\eps (t,x) +i \frac{\cos
t}{R}\nabla\chi\left(\frac{x}{R}\right) 
u^\eps(t,x),$$ 
and for $\frac{\pi}{2} -\Lambda\eps\leq t 
\leq \frac{\pi}{2} +\Lambda\eps$,
$$\left\| \frac{\cos
t}{R}\nabla\chi\left(\frac{\cdot}{R}\right) 
u^\eps(t,.)\right\|_{L^2}\leq C \frac{|\cos t|}{R}\leq
C\frac{\Lambda\eps}{R}.$$ 
Therefore to prove Lemma~\ref{lem:troncature} when $A^\eps =J^\eps$,
it is enough to prove,
 $$\limsup_{\eps \rightarrow 0} \sup_{\frac{\pi}{2} -\Lambda\eps\leq t 
\leq \frac{\pi}{2} +\Lambda\eps}  \left\|\left(1-\chi \left(\frac{\cdot}{R}\right)
\right)J^\eps(t) u^\eps (t,.)
\right\|_{L^2}= 0.$$
The function $J^\eps(t) u^\eps$ satisfies, from the commutation
property (\ref{eq:commut}),
\begin{equation}\label{decadix}
\left(i\eps \d_t +\frac{1}{2}\eps^2\Delta -\frac{|x|^2}{2} \right)J^\eps(t)
u^\eps = \eps^\sigma J^\eps(t) \left( |u^\eps|^{2\sigma}
 u^\eps \right).
\end{equation}
Notice that from Prop.~\ref{prop:avt} and
(\ref{eq:bornesnl}), Sobolev inequality implies that there exists
$C=C(\Lambda)$
such that for any $t \in [0, \pi/2 +\Lambda \eps]$,
\begin{equation}\label{eq:Ls}
\| u^\eps (t)\|_{L^{\underline s}} \leq \frac{C}{(|\cos t|
+\eps)^{\delta({\underline s})}}.
\end{equation}
At this stage, $C$ might depend on $\Lambda$ (even though we will know
it does not, afterward). 
Therefore, Cor.~\ref{cor:estgen}, 
applied to (\ref{decadix}) a finite number of times to 
cover the interval $[\frac{\pi}{2}-\Lambda_0
\eps,\frac{\pi}{2}+\Lambda \eps]$, implies that for any
$\Lambda \geq \Lambda_0$, $J^\eps(t)
u^\eps$ is bounded in $L^2$ for $t \in [0, \pi/2 +\Lambda
\eps]$. Next, commuting the cut-off function $\chi$ with
(\ref{decadix}) yields, 
\begin{equation*}
\begin{aligned}
\left(i\eps \d_t +\frac{1}{2}\eps^2\Delta -\frac{|x|^2}{2} \right)&
\left(1-\chi \left(\frac{x}{R}\right)\right)J^\eps(t)
u^\eps =  \eps^{n\sigma}\left(1-\chi \left(\frac{x}{R}\right)\right) J^\eps(t)
\left( |u^\eps|^{2\sigma} 
 u^\eps \right)\\
& -\frac{\eps^2}{2R^2}\Delta\chi \left(\frac{x}{R}\right)J^\eps(t)
u^\eps -\frac{\eps^2}{R}\nabla\chi \left(\frac{x}{R}\right)\nabla J^\eps(t)
u^\eps.
\end{aligned}
\end{equation*}
From Cor.~\ref{cor:estgen}  and (\ref{eq:Ls}), if we denote
$t^\eps_-=\pi/2-\Lambda\eps$, $t^\eps_+=\pi/2+\Lambda \eps$, we have
\begin{equation*}
\begin{aligned}
\left\| \left(1-\chi \left(\frac{.}{R}\right)\right)J^\eps(t)
u^\eps\right\|_{L^\infty(t^\eps_-,t^\eps_+;L^2)} & \leq
\left\| \left(1-\chi \left(\frac{\cdot}{R}\right)\right)J^\eps(t^\eps_-)
u^\eps\right\|_{L^2}\\ 
 +C\frac{\eps}{2R^2} \|J^\eps(t)
u^\eps\|_{L^1(t^\eps_-,t^\eps_+;L^2)}
& + C \frac{\eps}{R}\left\|\nabla J^\eps(t)u^\eps
\right\|_{L^1(t^\eps_-,t^\eps_+;L^2)}\\
& \leq \left\| \left(1-\chi \left(\frac{\cdot}{R}\right)\right)J^\eps(t^\eps_-)
u^\eps\right\|_{L^2}\\ 
+ C\Lambda\frac{\eps^2}{2R^2} \|J^\eps(t)
u^\eps\|_{L^\infty(t^\eps_-,t^\eps_+;L^2)}
& + C\Lambda \frac{\eps^2}{R}\left\|\nabla J^\eps(t)u^\eps
\right\|_{L^\infty(t^\eps_-,t^\eps_+;L^2)}.
\end{aligned}
\end{equation*}
We can conclude with the following lemma, whose proof is postponed to
Sect.~\ref{sec:ordre2}. 
\begin{lem}\label{lem:mulder}
Assume $f\in {\EuScript H}$ and $\sigma >1/2$. 
Let $\Lambda >1$. There exists $C=C(\Lambda)$ such that for any $t \in
[\pi/2-\Lambda \eps, \pi/2 +\Lambda \eps]$,
\begin{equation*}
\left\|\eps\nabla J^\eps(t)u^\eps
\right\|_{L^2} + \left\|\eps\nabla H^\eps(t)u^\eps
\right\|_{L^2} \leq C.
\end{equation*}
\end{lem}
This completes the proof of Lemma~\ref{lem:troncature}, the
computations with $H^\eps$ being similar.  
\end{proof}

\noindent To prove Prop.~\ref{prop:traversee1}, we now have to compare
$\tilde v^\eps$ and the truncated exact solution $u_R^\eps$. 

\begin{lem}\label{lem:mouse}
Assume $f\in {\EuScript H}$ and 
take $R=\eps^\alpha$. Then for any $0<\alpha <1$,
\begin{equation*}
\limsup_{\eps \rightarrow 0} \sup_{\frac{\pi}{2} -\Lambda\eps\leq t 
\leq \frac{\pi}{2} +\Lambda\eps}  \|A^\eps(t)(u_R^\eps(t)-\tilde v^\eps(t)
)\|_{L^2}\Tend \Lambda {+\infty} 0,
\end{equation*}
where $A^\eps(t)$ is either of the operators $\Id$, $J^\eps(t)$ or
$H^\eps(t)$.
\end{lem}
\begin{proof}[Proof of Lemma~\ref{lem:mouse}]
 Denote $\tilde
w^\eps_R=u_R^\eps-\tilde v^\eps$. Recall that $u_R^\eps$ solves 
\begin{equation*}
\left(i\eps \d_t +\frac{1}{2}\eps^2\Delta -\frac{|x|^2}{2} \right)u_R^\eps=
\eps^{n\sigma} |u^\eps|^{2\sigma}u_R^\eps+\frac{\eps^2}{2R}\nabla\chi\left(
\frac{x}{R}\right)\cdot\nabla u^\eps +\left( \frac{\eps}{R}\right)^2 
\Delta\chi\left( \frac{x}{R}\right) u^\eps,
\end{equation*}
and notice that with our choice for the cut-off function $\chi$,
$$\chi\left(\frac{x}{R}\right)=\chi\left(\frac{x}{2R}\right)
\chi\left(\frac{x}{R}\right),$$
therefore
\begin{equation*}
\begin{aligned}
\left(i\eps \d_t
+\frac{1}{2}\eps^2\Delta\right)u_R^\eps=&V_R(x)u_R^\eps+
\eps^{n\sigma} |u^\eps|^{2\sigma}u_R^\eps\\
 &+\frac{\eps^2}{2R}\nabla\chi\left(
\frac{x}{R}\right)\cdot\nabla u^\eps +\left( \frac{\eps}{R}\right)^2 
\Delta\chi\left( \frac{x}{R}\right) u^\eps,
\end{aligned}
\end{equation*}
where 
$$V_R(x) =\chi\left(\frac{x}{2R}\right) 
\frac{|x|^2}{2}.$$
The remainder $\tilde w^\eps_R$ solves
\begin{equation}\label{eq:resteb}
\begin{aligned}
\left(i\eps \d_t
+\frac{1}{2}\eps^2\Delta\right)\tilde w_R^\eps=&V_R(x)u_R^\eps+
\eps^{n\sigma} \left(|u^\eps|^{2\sigma}u_R^\eps - |\tilde
v^\eps|^{2\sigma}\tilde v^\eps\right)\\ 
&+\frac{\eps^2}{2R}\nabla\chi\left(
\frac{x}{R}\right)\cdot\nabla u^\eps +\left( \frac{\eps}{R}\right)^2 
\Delta\chi\left(\frac{x}{R}\right) u^\eps.
\end{aligned}
\end{equation}
Apply Prop.~\ref{prop:estgen}, with now ${\EuScript U}^\eps
=U_0^\eps$, $F^\eps=0$ and 
$h^\eps =h^\eps_1+h^\eps_2$, where
$$h^\eps_1 = V_R(x)u_R^\eps  +\frac{\eps^2}{2R}\nabla\chi\left(
\frac{x}{R}\right)\cdot\nabla u^\eps +\left( \frac{\eps}{R}\right)^2 
\Delta\chi\left(\frac{x}{R}\right) u^\eps,$$
and 
$$h^\eps_2=\eps^{n\sigma}\left(|u^\eps|^{2\sigma}u_R^\eps - |\tilde
v^\eps|^{2\sigma}\tilde v^\eps\right).$$
This yields, for $\pi/2-\Lambda \eps\leq t_0\leq t_1\leq \pi/2+\Lambda
\eps$, 
\begin{equation*}
\begin{aligned}
\left\|\tilde w_R^\eps \right\|_{L^{\underline
q}(t_0,t_1;L^{\underline r})}&  \leq C\eps^{-\frac{1}{\underline
q}}\left\|\tilde w_R^\eps(t_0)\right\|_{L^2}+
C\eps^{-1-\frac{1}{\underline q}}\|h^\eps_1\|_{L^1(t_0,t_1;L^2)}\\
& + C\eps^{-1-\frac{2}{\underline q}}\|h^\eps_2\|_{L^{{\underline q}'}
(t_0,t_1;L^{{\underline r}'})}\\
& \leq C\eps^{-\frac{1}{\underline
q}}\left\|\tilde w_R^\eps(t_0)\right\|_{L^2}+
C\eps^{-1-\frac{1}{\underline q}}
\int_{t_0}^{t_1}\left(R^2+\frac{\eps}{R}+
\left(\frac{\eps}{R}\right)^2  \right)dt\\
& +{\underline C}\left(
\frac{t_1-t_0}{\eps}\right)^{2\sigma/{\underline k}}
\left(\left\|\tilde w_R^\eps
\right\|_{L^{\underline q}(t_0,t_1;L^{\underline r})}+\left\| w_R^\eps
\right\|_{L^{\underline q}(t_0,t_1;L^{\underline r})} \right),
\end{aligned}
\end{equation*}
from H\"older inequality. Taking 
$${\underline C}\left(
\frac{t_1-t_0}{\eps}\right)^{2\sigma/{\underline k}}\leq
\frac{1}{2},$$
we have 
\begin{equation*}
\begin{aligned}
\left\|\tilde w_R^\eps \right\|_{L^{\underline
q}(t_0,t_1;L^{\underline r})}&  \leq C\eps^{-\frac{1}{\underline
q}}\left\|\tilde w_R^\eps(t_0)\right\|_{L^2}+
C\eps^{-\frac{1}{\underline q}}
\left(R^2+\frac{\eps}{R}+
\left(\frac{\eps}{R}\right)^2  \right)\\
& +C\left\| w_R^\eps
\right\|_{L^{\underline q}(t_0,t_1;L^{\underline r})}.
\end{aligned}
\end{equation*}
Repeating this manipulation a finite number of times covers the whole
interval $t\in [\pi/2-\Lambda_0\eps,\pi/2+\Lambda_0 \eps]$. Doing
this, we get a possibly large, but finite, constant, which can be seen
as the analogue of the exponential term in Gronwall lemma. For
$\Lambda_0 \eps\leq |t-\pi/2|\leq \Lambda \eps$, we 
use time decay estimates provided by
$J^\eps$; bearing the comparison with Gronwall lemma in mind, this
means that the operator $J^\eps$ provides some integrability for
$\Lambda_0 \eps\leq |t-\pi/2|$, which is stated in (\ref{eq:hypF}),
and implies the first condition in Cor.~\ref{cor:estgen}. This
integrability is needed to get a 
bound independent of $\Lambda \geq \Lambda_0$.    
When $A^\eps=H^\eps$, from (\ref{eq:commut}),
$$\left[H^\eps(t), i\eps \d_t +\frac{1}{2}\eps^2\Delta 
\right]=\left[H^\eps(t), \frac{|x|^2}{2}
\right]= i\eps x\sin t,$$
and $ H^\eps u_R^\eps$ satisfies,
\begin{equation*}
\begin{aligned}
\left(i\eps \d_t +\frac{1}{2}\eps^2\Delta \right)H^\eps
u_R^\eps= & i\eps x(\sin t) u_R^\eps +
\eps^{n\sigma} H^\eps\left( |u^\eps|^{2\sigma}u_R^\eps\right)\\
&+ V_R(x)H^\eps u_R^\eps + i\eps (\sin t)\nabla V_R (x)u_R^\eps\\
& +\frac{\eps^2}{2R}\nabla\chi\left(
\frac{x}{R}\right)H^\eps\d_x u^\eps +i\eps\sin t \frac{\eps^2}{2R^2}
\Delta\chi\left(
\frac{x}{R}\right)\nabla u^\eps \\
& +\left( \frac{\eps}{R}\right)^2 \Delta\chi\left(
\frac{x}{R}\right) H^\eps u^\eps+i\eps\sin t
 \frac{\eps^2}{R^3}\nabla \Delta \chi\left(
\frac{x}{R}\right) u^\eps.
\end{aligned}
\end{equation*}
It follows that the remainder $H^\eps \tilde w_R^\eps$ satisfies,
\begin{equation*}
\begin{aligned}
\left(i\eps \d_t +\frac{1}{2}\eps^2\Delta \right)H^\eps
\tilde w_R^\eps= & 
\eps^{n\sigma} H^\eps\left( |u^\eps|^{2\sigma}u_R^\eps - |\tilde
v^\eps|^{2\sigma}\tilde v^\eps\right)\\
& +ix\eps\chi\left(\frac{x}{R}\right)(\sin t) u^\eps+i\eps (\sin
t)\nabla V_R (x)u_R^\eps \\
&+ V_R(x) \chi\left(
\frac{x}{R}\right)H^\eps u^\eps+i \frac{\eps}{R}V_R(x)(\sin t)\nabla\chi\left(
\frac{x}{R}\right)u^\eps\\
& +\frac{\eps^2}{2R}\nabla\chi\left(
\frac{x}{R}\right)H^\eps\nabla u^\eps +i\eps\sin t \frac{\eps^2}{2R^2}
\Delta\chi\left(
\frac{x}{R}\right)\nabla u^\eps \\
& +\left( \frac{\eps}{R}\right)^2 \Delta\chi\left(
\frac{x}{R}\right) H^\eps u^\eps+i\eps\sin t
 \frac{\eps^2}{R^3} \nabla\Delta\chi\left(
\frac{x}{R}\right) u^\eps.
\end{aligned}
\end{equation*}
We can estimate the term in $H^\eps\d_x u^\eps$ because we can
estimate $\d_x H^\eps u^\eps$ (Lemma~\ref{lem:mulder}) and the
following holds, 
$$[H^\eps(t),\nabla]=-\cos t=O(\eps) \textrm{ for }\pi/2 -\Lambda \eps
\leq t \leq \pi/2 +\Lambda \eps .$$
The proof then proceeds as above. 
\end{proof}

Lemmas~\ref{lem:troncature} and \ref{lem:mouse} clearly imply
Prop.~\ref{prop:traversee1}.
\end{proof}

The assumption $f\in {\EuScript
H}$ turns out to be unnecessary. Indeed, we can use a density argument
for $\tilde v^\eps$, and approach $f \in \Sigma$ by functions in ${\EuScript
H}$ up to a small error in the norms that are considered in
Prop.~\ref{prop:traversee1}; this stems from
global well-posedness of (\ref{eq:nls}) (see e.g. \cite{Caz},
\cite{Ginibre}).  We can mimic the proof of this result for $u^\eps$,
thanks to $J^\eps$ and $H^\eps$. 

\begin{prop}\label{prop:traversee}
Proposition~\ref{prop:traversee1} still holds if we
assume $f\in \Sigma$ and $r^\eps\not \equiv 0$. 
\end{prop}

\subsection{Proof of Lemma~\ref{lem:mulder}}
\label{sec:ordre2}

We first use the following remark. 
\begin{lem}
Assume a function $u^\eps(x)$ satisfies
$$\|(-\eps^2\Delta +|x|^2)u^\eps\|_{L^2}+ \|\eps\nabla u^\eps\|_{L^2} + 
\|xu^\eps\|_{L^2}\leq \underline{C},$$
where $\underline C$ does not depend on $\eps$. Then 
$$\|\eps^2\Delta u^\eps\|_{L^2}+\||x|^2 u^\eps\|_{L^2}\leq \underline C.$$
\end{lem}
Now the idea is to differentiate (\ref{eq:nl}) with respect to
time. This is classical for the case of the nonlinear Schr\"odinger
equation (\ref{eq:nls}), see e.g. \cite{Caz}, Sect.~5.2. Thanks to the
above lemma, we can adapt the mentioned results to prove the
following proposition.
\begin{prop}\label{prop:H2}
Assume $f\in {\EuScript H}$. Let $\Lambda >1$. Then $$u^\eps \in
C\left(0,\frac{\pi}{2}+\Lambda\eps;H^2\cap {\F}(H^2)\right)\cap
C^1\left(0,\frac{\pi}{2}+\Lambda\eps;L^2\right),$$ 
and there exists 
$C=C(\Lambda)$ independent of $\eps$ such that
$$\sup_{0\leq t\leq \frac{\pi}{2}+\Lambda\eps}\| \eps\d_t u^\eps(t)
\|_{L^2}+\sup_{0\leq t\leq \frac{\pi}{2}+\Lambda\eps}\| \eps^2\Delta
u^\eps(t) 
\|_{L^2}+\sup_{0\leq t\leq \frac{\pi}{2}+\Lambda\eps}\| |x|^2 u^\eps(t)
\|_{L^2} \leq C.$$
\end{prop}
\begin{proof}[Idea of the proof]
As in \cite{Caz}, the idea of the
proof consists in differentiating the equation satisfied by $u^\eps$
with respect to time, and estimate $\eps\d_t u^\eps$. Since the
harmonic potential commutes with the time derivative, one can mimic
the proof given in \cite{Caz}, Sect.~5.2. 
When there is no potential, like in (\ref{eq:nls}), the control of the
nonlinear term and the time derivative give some control on
$\eps^2\Delta u^\eps$. In our case, this controls $\eps^2\Delta u^\eps
-|x|^2 u^\eps$. From the above lemma, this means that we can estimate
each of these two terms. \qed

Notice that the following algebraic
identity holds point-wise, for any $j,k$,
\begin{equation}\label{eq:alg2}
\begin{aligned}
&|x_k\eps J_j^\eps(t)u^\eps|^2 +|\eps^2\d_k J_j^\eps(t)u^\eps|^2 +
 |x_kH_j^\eps(t)u^\eps|^2 +|\eps \d_k H_j^\eps(t)u^\eps|^2=\\
&= |x_jx_k u^\eps(t)|^2 +|\eps x_k\d_j u^\eps(t)|^2+ |\eps
\delta_{jk}u^\eps(t)+\eps x_j\d_k u^\eps(t)|^2+|\eps^2\d_{j
k}^2u^\eps(t)|^2, 
\end{aligned}
\end{equation}
where $\delta_{jk}$ stands for the Kronecker symbol. From
Prop.~\ref{prop:H2}, the right-hand side is bounded in $L^1_x$,
uniformly for $0\leq t\leq \frac{\pi}{2}+\Lambda\eps$. This implies
the boundedness of $\eps\nabla H^\eps(t)u^\eps$ stated in
Lemma~\ref{lem:mulder}, and even a little more, that is,
\begin{equation}\label{marre}
\forall t\in \left[0,\frac{\pi}{2}+\Lambda \eps\right],\ \
\|\eps\nabla H^\eps(t)u^\eps\|_{L^2}\leq C. 
\end{equation}
At this stage, we have not assumed that the nonlinearity was twice
differentiable. 
On the other hand, we just have
$$\|\eps\nabla J^\eps(t)u^\eps\|_{L^2}\leq \frac{C}{\eps}.$$
The idea is that for this term, (\ref{eq:alg2}) is far from giving
a sharp estimate. Indeed, for $|t-\pi/2|=O(\eps)$, we guess that the
main contribution of $u^\eps$ lies in $|x|=O(\eps)$ (semi-classical
Schr\"odinger equations are morally hyperbolic). This is precisely what
we have to prove. With the additional
remark that near $t=\pi/2$, one can replace $H^\eps$ with $\eps\nabla$
up to a small error term, this suggests that the leading order term of
the left-hand side is $|\eps \d_{x_k} H_j^\eps(t)u^\eps|^2$, and the 
leading order term of
the right-hand side is $|\eps^2\d_{x_j x_k}^2u^\eps(t)|^2$. Thus there is
nothing more to hope from this identity. 

This in fact must not be surprising. The only additional estimates we
obtained are those stated in Prop.~\ref{prop:H2},
$$\sup_{0\leq t\leq \frac{\pi}{2}+\Lambda\eps}\| \eps\d_t u^\eps(t)
\|_{L^2}+\sup_{0\leq t\leq \frac{\pi}{2}+\Lambda\eps}\| \eps^2\Delta
u^\eps(t) 
\|_{L^2}+\sup_{0\leq t\leq \frac{\pi}{2}+\Lambda\eps}\| |x|^2 u^\eps(t)
\|_{L^2} \leq C.$$
The boundedness of the first two terms means that $u^\eps$ is 
$\eps$-oscillating, and the boundedness of the last term means that
the solution remains confined. This is due to the fact that we work
with an unbounded potential, but not to the fact that we consider the
harmonic potential in particular. Therefore, there is no precise
geometric information in this estimate. As a matter of fact, away from
the focus, this kind
of information is given by the operator $J^\eps$.
  
We assume that
the nonlinearity is twice differentiable. 
Recall that for $1\leq j\leq n$, $J^\eps_j(t)u^\eps$ satisfies
$$\left(i\eps \d_t +\frac{1}{2}\eps^2\Delta -\frac{|x|^2}{2}
\right)J_j^\eps(t)u^\eps = \eps^{n\sigma}J_j^\eps(t)\left(|u^\eps|^{2\sigma}
u^\eps\right).$$
Differentiating this equation with respect to $x_k$ yields
\begin{equation}\label{eq:rew}
\left(i\eps \d_t +\frac{1}{2}\eps^2\Delta -\frac{|x|^2}{2}
\right)\eps\d_kJ_j^\eps(t)u^\eps =
\eps^{n\sigma+1}\d_kJ_j^\eps(t)\left(|u^\eps|^{2\sigma} 
u^\eps\right)+\eps x_k J_j^\eps(t)u^\eps. 
\end{equation}
The last term comes from the commutation of the harmonic potential
with $\d_{x_k}$. From (\ref{eq:jauge}), the following point-wise
estimate holds,
\begin{equation}
\left|\eps \d_kJ_j^\eps(t)\left(|u^\eps|^{2\sigma} 
u^\eps\right)\right| \lesssim |u^\eps|^{2\sigma-1} |\eps\d_k
u^\eps|. |J_j^\eps(t)u^\eps|+ |u^\eps|^{2\sigma}
|\eps\d_kJ_j^\eps(t)u^\eps|.  
\end{equation}
The idea is that the last term is well prepared to apply Gronwall
lemma. For the first term of the left-hand side, we have to work a
little more. Apply Cor.~\ref{cor:estgen} to (\ref{eq:rew}), with
\begin{equation*}
\begin{aligned}
&|F^\eps(\eps\d_kJ_j^\eps(t)u^\eps)|\lesssim |u^\eps|^{2\sigma}
|\eps\d_kJ_j^\eps(t)u^\eps|,\\
&|h^\eps|\lesssim
\eps^{n\sigma}|u^\eps|^{2\sigma-1} 
|\eps\d_k u^\eps|. |J_j^\eps(t)u^\eps|+|\eps x_k J_j^\eps(t)u^\eps|.
\end{aligned}
\end{equation*}
We already know that $F^\eps$ satisfies
(\ref{eq:hypF}). Cor.~\ref{cor:estgen} yields,
\begin{equation}\label{eq:16h35}
\begin{aligned}
\left\|\eps\d_kJ_j^\eps(t)u^\eps\right\|_{L^\infty(t_0,t_1;L^2)}& \leq 
\|\eps \d_k(x_j f)\|_{L^2}  +C\left\|x_k
J_j^\eps(t)u^\eps\right\|_{L^1(t_0,t_1;L^2)}\\
&+C\eps^{n\sigma -1-\frac{1}{\underline q}} 
\left\||u^\eps|^{2\sigma-1} 
\eps\d_k u^\eps J_j^\eps(t)u^\eps \right\|_{L^{{\underline
q}'}(t_0,t_1;L^{{\underline r}'})} .
\end{aligned}
\end{equation}
For fixed $t$, H\"older inequality yields,
$$\left\||u^\eps|^{2\sigma-1} 
\eps\d_k u^\eps J_j^\eps(t)u^\eps \right\|_{L^{{\underline r}'}}\leq
\| u^\eps\|_{L^{a_1}}^{2\sigma-1}\|\eps\d_k u^\eps\|_{L^{a_2}}\|J_j^\eps(t)
u^\eps\|_{L^{\underline r}}\ \ ,$$
with
$$\frac{2\sigma-1}{a_1}+\frac{1}{a_2}=\frac{2\sigma}{\underline s}.$$
We can take for instance $a_1=a_2=\underline s$.
This implies, along with (\ref{eq:sobolevn}), since $\|
\eps\d_k H^\eps(t)u^\eps\|_{L^2}$ is uniformly bounded 
$$\left\||u^\eps|^{2\sigma-1} 
\eps\d_k u^\eps J_j^\eps(t)u^\eps \right\|_{L^{{\underline r}'}}\leq
\frac{C}{(\cos t +\eps)^{2\sigma\delta (\underline s)}}\left(\|
J^\eps(t)\eps\d_k u^\eps\|_{L^2}^{\delta (a_2)}+ 1\right)\|J_j^\eps(t)
u^\eps\|_{L^{\underline r}}.$$
Now apply H\"older inequality in time, with
$$\frac{1}{\underline q '}=\frac{2\sigma}{\underline k}+
\frac{1}{\infty}+\frac{1}{\underline q}.$$
This yields
\begin{equation}
\begin{aligned}
\bigl\||u^\eps|^{2\sigma-1} 
\eps\d_k u^\eps &J_j^\eps(t)u^\eps \bigr\|_{L^{{\underline
q}'}(t_0,t_1;L^{{\underline r}'})} \lesssim   \ A^\eps (t_0,t_1)
\|J_j^\eps(t) 
u^\eps\|_{L^{\underline q}(t_0,t_1;L^{\underline r})}\\
+ & A^\eps (t_0,t_1)\|
J^\eps(t)\eps\d_k u^\eps\|_{L^\infty(t_0,t_1;L^2)}^{\delta (a_2)}
\|J_j^\eps(t) 
u^\eps\|_{L^{\underline q}(t_0,t_1;L^{\underline r})},
\end{aligned}
\end{equation}
where $A^\eps$ is defined in Prop.~\ref{prop:estgen}. We also know
that
$$\|J_j^\eps(t) 
u^\eps\|_{L^{\underline q}(t_0,t_1;L^{\underline r})} \leq C
\eps^{-1/\underline q},$$
therefore (\ref{eq:16h35}) yields,
\begin{equation*}
\begin{aligned}
\bigl\|\eps\d_kJ_j^\eps(t)u^\eps&\bigr\|_{L^\infty(t_0,t_1;L^2)}\leq  
\|\eps \d^2_{jk} f\|_{L^2}  +C \left\|x_k
J_j^\eps(t)u^\eps\right\|_{L^1(t_0,t_1;L^2)}\\
&+C\eps^{n\sigma -1-\frac{2}{\underline q}} A^\eps (t_0,t_1) \left(1 + 
\|J^\eps(t)\eps\d_k u^\eps\|_{L^\infty
(t_0,t_1;L^2)}^{\delta (a_2)}\right).
\end{aligned}
\end{equation*}
But from Lemma~\ref{lem:alg},
$$n\sigma -1-\frac{2}{\underline q}= 2\sigma \left(\delta(\underline
s) - \frac{1}{\underline k} \right),$$
and we find the same quantity as in Prop.~\ref{prop:estgen}, that is
$\eps^{2\sigma \left(\delta(\underline
s) - \frac{1}{\underline k} \right)}A^\eps (t_0,t_1).$ With the remarks
that 
$$\left[\eps\d_k, J_j^\eps(t) \right]= \delta_{jk}\sin t,$$
and $\delta(a_2)\leq 1$, 
we have also,
\begin{equation}\label{eq:ven13}
\begin{aligned}
\bigl\|\eps\d_kJ_j^\eps(t)u^\eps&\bigr\|_{L^\infty(t_0,t_1;L^2)}\leq  
\|\eps \d^2_{jk} f\|_{L^2}  +C \left\|x_k
J_j^\eps(t)u^\eps\right\|_{L^1(t_0,t_1;L^2)}\\
&+C\eps^{n\sigma -1-\frac{2}{\underline q}} A^\eps (t_0,t_1) \left(1 + 
\|\eps\d_k J^\eps(t)u^\eps\|_{L^\infty
(t_0,t_1;L^2)}\right).
\end{aligned}
\end{equation}
Now it is natural to study $x_k J_j^\eps(t)u^\eps$. It satisfies,
\begin{equation*}
\left(i\eps \d_t +\frac{1}{2}\eps^2\Delta -\frac{|x|^2}{2}
\right)x_kJ_j^\eps(t)u^\eps =
\eps^{n\sigma}x_kJ_j^\eps(t)\left(|u^\eps|^{2\sigma} 
u^\eps\right)+\eps^2 \d_k J_j^\eps(t)u^\eps. 
\end{equation*}
The same computation as above, minus the three terms estimate which
is not needed here, yields
\begin{equation}\label{eq:ven13'}
\begin{aligned}
\bigl\|x_kJ_j^\eps(t)u^\eps&\bigr\|_{L^\infty(t_0,t_1;L^2)}\leq  
\|x_k \d_j f\|_{L^2}  +C \left\|\eps \d_k
J_j^\eps(t)u^\eps\right\|_{L^1(t_0,t_1;L^2)}\\
&+C\eps^{n\sigma -1-\frac{2}{\underline q}} A^\eps (t_0,t_1)  
\|x_k J_j^\eps(t)u^\eps\|_{L^\infty
(t_0,t_1;L^2)}.
\end{aligned}
\end{equation}
Summing (\ref{eq:ven13}) and (\ref{eq:ven13'})
over $j$ and $k$, Lemma~\ref{lem:mulder} follows from the
Gronwall lemma.
\end{proof}

\subsection{Past the first focus}
\label{sec:past}
After the first focus, we can proceed like before the focus, and
iterate this process. Notice that if $n\geq 3$, then
$\frac{1}{2}>\sigma_0(n)$, and we always have $\sigma >\sigma_0(n)$.
Next, we can prove the analogous of
Prop.~\ref{prop:avt}, using Prop.~\ref{prop:matching}
and Corollary~\ref{cor:idem}.
\begin{prop}\label{prop:apres}
The following asymptotics holds for $\pi/2< t \leq \pi$, 
\begin{equation*}
\limsup_{\eps \rightarrow 0}
\sup_{\pi/2+\Lambda \eps \leq t \leq \pi}\|A^\eps(t)(u^\eps
-v_1^\eps)(t)\|_{L^2}
\Tend \Lambda {+\infty} 0,
\end{equation*}
where $A^\eps(t)$ is either of the operators $\Id$, $J^\eps(t)$ or
$H^\eps(t)$.
\end{prop}
Finally, $v_{{\rm app},1}^\eps$ approximates $v_1^\eps$ like in
Corollary~\ref{cor:idem}. Corollary~\ref{cor:idem}, 
Prop.~\ref{prop:avt},
\ref{prop:traversee} and \ref{prop:apres} then imply
Theorem~\ref{th:princ}.  

When $t=\pi$, the problem is almost the
same as at time $t=0$. The initial data $f$ is replaced by $f_1$, and 
$$\| u^\eps (\pi,.)- f_1\|_\Sigma =o(1).$$
Therefore, Theorem~\ref{th:princ} can be iterated, which
yields Corollary~\ref{cor:iter}, because of the property, 
$$\forall \theta \in \R,\ \forall \psi_-\in \Sigma, \ 
S(e^{i\theta}\psi_-)=e^{i\theta}S(\psi_-).$$

\section{When the nonlinearity is focusing}
\label{sec:blowup}
In this section, we assume $n=1$ for simplicity. 
The first remark to guess the result of Prop.~\ref{prop:blowup} is
that in the proof of Prop.~\ref{prop:avt}, the sign of the
nonlinearity in unimportant. One needs local existence results to start
the ``so long'' argument, and general estimates on the
nonlinear term that do not involve its sign. Therefore
Prop.~\ref{prop:avt} still holds when $u^\eps$ is the solution of
(\ref{eq:nlfoc}). 

Next, assume for a moment that the matching argument can be used as in
Prop.~\ref{prop:matching}, and that afterward, the harmonic potential
can be neglected because of concentration. The behavior of $u^\eps$
should then be the same as the solution of 
$$i\eps\d_t v^\eps +\frac{1}{2}\eps^2\d_x^2 v^\eps 
=-\eps^2 |v^\eps|^4 v^\eps.$$
Resuming the scaling (\ref{eq:scaling}), we have to understand the
behavior of the solution of the same equation with $\eps =1$. It is
well known (see \cite{Caz}) that for small initial data, the solution
exists globally. The critical mass is the $L^2$-norm of the ground state
$R$ defined in Prop.~\ref{prop:blowup}. Recall what happens in this
critical case.
\begin{theo}(\cite{MerleDuke}, case $n=1$) Let $\varphi \in \Sigma$, with
$\|\varphi \|_{L^2}=\|R \|_{L^2}$. Let $\psi$ be the solution of the
initial value problem,
\begin{equation*}
\left\{
\begin{aligned}
i\d_t \psi +\frac{1}{2}\d_x^2 \psi &
= - |\psi|^4 \psi,\\ 
\psi_{\mid t=0} & = \varphi.
\end{aligned}
\right.
\end{equation*}
Assume that $\psi$ blows up at time $t=t_*$. Then there exist
$\theta$, $\omega$, $\xi_0$, $x_1\in \R$ such that for $t<t_*$,
\begin{equation}\label{eq:profmerle}
\psi(t,x)= \sqrt{\frac{\omega}{t_*-t}}R\left(\omega \left(
\frac{x-x_1}{t_*-t}-\xi_0 \right) \right)e^{i\left(\theta+ 
\frac{\omega^2}{t_*-t}-\frac{(x-x_1)^2}{2(t_*-t)} \right)}.
\end{equation}
\end{theo}
The second important remark is that such profiles as in
(\ref{eq:profmerle}) are dispersed when $t$ goes to $-\infty$. If
$\omega =1$, $x_1=\xi_0=0$, then 
\begin{equation}\label{eq:&}
\left\| U_0(t_*-t)\psi (t) -\frac{1}{\sqrt{2i\pi}}\widehat R
\right\|_\Sigma \Tend t {-\infty} 0.
\end{equation}
From the uniqueness in the
first part of Prop.~\ref{prop:scattnls}, if $\psi$ solves the critical
nonlinear Schr\"odinger equation and behaves asymptotically when
$t\rightarrow -\infty$ like the free evolution of
$$\psi_-^{t_*}:= \frac{1}{\sqrt{2i\pi}}U_0(-t_*)\widehat R,$$
then $\psi$ is given by (\ref{eq:profmerle}) with $\omega =1$ and
$x_1=\xi_0=0$. Back to the scaling (\ref{eq:scaling}), this yields the
definitions $f(x)= R(x)e^{i\frac{t_*}{2}|x|^2}$ (from (\ref{eq:psi-}) and the
definition of $\psi_-^{t_*}$) and  (\ref{eq:solappfoc}).

Now sketch the proof of Prop.~\ref{prop:blowup}. As we noticed,
Prop.~\ref{prop:avt} describes the behavior of $u^\eps$ up to
$t=\pi/2 -\Lambda\eps$ for large $\Lambda$. What prevents us from
mimicking the proof of Prop.~\ref{prop:matching}? The limit
(\ref{eq:step1}) still holds, as well as Lemmas~\ref{lem:1} and
\ref{lem:3}. However, one cannot apply Lemma~\ref{lem:2} so easily to
$u^\eps$ and $v^\eps_{\rm app}$ for estimate (\ref{eq:bornesnl}) is
not true when the nonlinearity is focusing. On the other hand,
(\ref{eq:bornesnl}) is true up to time  $t=\pi/2 -\Lambda\eps$ for
large $\Lambda$, from  Prop.~\ref{prop:avt} and the algebraic
identity (\ref{eq:algebre}). Therefore Prop.~\ref{prop:matching} still
holds. 

Finally, one can adapt Prop.~\ref{prop:traversee1} by replacing the
time interval $[\pi/2 -\Lambda \eps, \pi/2 +\Lambda \eps]$ by 
$[\pi/2 -\Lambda \eps, \pi/2 +t_*\eps-\lambda\eps]$, for any positive
$\lambda$. The method of our proof does not allow to go
further. Indeed, we have the following estimates,
$$\left\| \tilde v^\eps \left(\frac{\pi}{2} +t_*\eps-\lambda\eps \right)
\right\|_{L^\infty} =\frac{\|R\|_{L^\infty}}{\sqrt{\lambda\eps}}, \ \ 
\left\| \eps\d_x \tilde v^\eps \left(\frac{\pi}{2} +t_*\eps-\lambda\eps \right)
\right\|_{L^2} =\frac{C(R)}{\lambda}.$$
Therefore, one cannot hope that (\ref{eq:solong2}) holds beyond
$t=\frac{\pi}{2} +t_*\eps-\lambda\eps$ (with $C_1$ proportional to
$\lambda^{-1/2}$). On the other hand, if our final time is $t=\frac{\pi}{2}
+t_*\eps-\lambda\eps$ with $\lambda >0$, we can prove the analogue of
Prop.~\ref{prop:traversee1} by a ``so long'' argument (that is,
(\ref{eq:solong2}) with $C_1$ proportional to
$\lambda^{-1/2}$). As a result, we have the first part of
Prop.~\ref{prop:blowup}. The last part follows from the remark we made
above, that we know $\tilde v^\eps$ explicitly, therefore in
particular its value at time $t=\frac{\pi}{2} +t_*\eps-\lambda\eps$. 
\section{Anisotropic harmonic potential}
\label{sec:anis}
Consider the general harmonic potential in $\R^n$,
\begin{equation}
\label{eq:harmogen}
V(x)=\frac{1}{2}(\omega_1^2 x_1^2+ \omega_2^2 x^2_2 + \ldots
+\omega_n^2 x_n^2), 
\end{equation}
with $\omega_j >0$ for all $j$. 
It is isotropic when all the $\omega_j$'s are equal, anisotropic
otherwise. We suppose that the $\omega_j$'s take exactly $d$ distinct
values ($2\leq d\leq n$), 
and renaming the space variables if necessary, we can assume
that
$$0<\omega_1< \omega_2< \ldots<\omega_d.$$
We denote $i_j$ the multiplicity of $\omega_j$, $1\leq j\leq d$
 ($i_1+\ldots +i_d=n$). At least two possibilities occur, as for the
result one can hope for, corresponding either to Th.~\ref{th:princ} or
to Cor.~\ref{cor:transl}. In the former case, one would be interested
in the Cauchy problem
\begin{equation}\label{eq:anis1}
\left\{
\begin{aligned}
i\eps\d_t u^\eps +\frac{1}{2}\eps^2\Delta u^\eps &
=V(x)u^\eps +\eps^{k\sigma} |u^\eps|^{2\sigma} u^\eps,\ \ \ \
(t,x)\in 
\R_+\times \R^n, \\ 
u^\eps_{\mid t=0} & = f(x) +r^\eps(x),
\end{aligned}
\right.
\end{equation}
and in the latter, in
\begin{equation}\label{eq:anis2}
\left\{
\begin{aligned}
i\eps\d_t u^\eps +\frac{1}{2}\eps^2\Delta u^\eps &
=V(x)u^\eps +\eps^{n\sigma} |u^\eps|^{2\sigma} u^\eps,\ \ \ \
(t,x)\in 
\R_+\times \R^n, \\ 
u^\eps_{\mid t=0} & = \frac{1}{\eps^{n/2}}f\left(\frac{x}{\eps}
\right) +\frac{1}{\eps^{n/2}}r^\eps\left(\frac{x}{\eps}
\right),
\end{aligned}
\right.
\end{equation}
where $f, r^\eps \in \Sigma$ and $\|r^\eps\|_\Sigma \Tend \eps
0 0$. We briefly discuss Eq.~(\ref{eq:anis1}), and explain more
precisely what happens for Eq.~(\ref{eq:anis2}).

For (\ref{eq:anis1}), the same method as in Sect.~\ref{sec:bkw} leads
to the following phase, profile and operators,
\begin{equation}
\begin{aligned}
\varphi(t,x)&= -\frac{1}{2} \sum_{j=1}^n \omega_j x_j^2 \tan (\omega_j
t),\\
v_0(t,x)&= \left(\prod_{j=1}^n \frac{1}{\sqrt{\cos (\omega_j t)}}
\right)f\left(\frac{x_1}{\cos (\omega_1 t)},\ldots, \frac{x_n}{\cos
(\omega_n t)}  \right),\\
J_j^\eps(t)&= \frac{\omega_j x_j}{\eps}\sin (\omega_j t) - i \cos
(\omega_j t)\d_j,\\
H_j^\eps(t)&= \omega_j x_j \cos (\omega_j t) +i \eps \sin
(\omega_j t)\d_j.
\end{aligned}
\end{equation}
The first focusing occurs for $t=\frac{\pi}{2\omega_d}$; the solution
$u^\eps$ focuses on the $i_d$-dimensional vector space defined by
$$E_d= \{ x_j=0 ,\ \forall j \textrm{ such that }\omega_j=\omega_d
\}.$$
Therefore, the critical index for the nonlinear term to be relevant in
(\ref{eq:anis1}) would be $k=\operatorname{dim} E_d= i_d$. If $k>i_d$,
then the nonlinear term 
remains negligible up to time $t=\frac{\pi}{2\omega_d}$ and before the
next focusing, where the same
discussion is valid. If $\dis k> \max_{1\leq j \leq d}i_j$, then the
nonlinear term is everywhere negligible, provided that no simultaneous
focusings occur; indeed, the $\omega_j$ part of the harmonic potential
will cause focusing at times 
$$\frac{\pi}{2\omega_j} +\frac{\kappa \pi}{\omega_j}, \ \kappa\in
\Z.$$
Two (or more) distinct $\omega_j$'s can cause cumulated focusing if
they are rationally related. To simplify the discussion, we now assume
$n=2$ and that $\omega_1$ and $\omega_2$ are irrationally related.
In that case, $u^\eps$ focuses at time
$t=\frac{\pi}{2\omega_2}$ on the line $\{x_2 = 0\}$. If $k=1$, then
the nonlinear term becomes relevant near $\{
(t,x_2)=(\frac{\pi}{2\omega_2},0)\}$. The case of a focusing on a line
was treated in \cite{Ca3} without potential,
with an initial oscillation that forces such a geometry for the
caustic. With an anisotropic oscillator, the situation is technically
much harder to 
handle. In \cite{Ca3}, no oscillation was present in the other space
variable, and this variable could be considered as a parameter. In
the present case, oscillations are always present in both space
variables, so it is harder to measure the dependence of $u^\eps$ with
respect to
$x_1$ when it focuses on $\{x_2 = 0\}$. We leave out the discussion at
this stage. \\

On the other hand, it is possible to understand (and prove) what
happens for Eq.~(\ref{eq:anis2}). Because we altered the time origin,
the operators we now use write,
\begin{equation}
\begin{aligned}
J_j^\eps(t)&= \frac{\omega_j x_j}{\eps}\cos (\omega_j t) + i \sin
(\omega_j t)\d_j,\\
H_j^\eps(t)&= \omega_j x_j \sin (\omega_j t) -i \eps \cos
(\omega_j t)\d_j.
\end{aligned}
\end{equation}
We also have an explicit formula for the linear solution (the analogue
of Eq.~(\ref{eq:mehler})), which yields
in particular Strichartz estimates. The solution of 
\begin{equation*}
\left\{
\begin{aligned}
i\eps\d_t v^\eps +\frac{1}{2}\eps^2\Delta v^\eps &
=V(x)v^\eps, \\
v^\eps_{\mid t=0}& = f(x),
\end{aligned}
\right.
\end{equation*}
is given by
\begin{equation*}
v^\eps(t,x) = \prod_{j=1}^n \left(\frac{\omega_j}{2i\pi\eps
\sin\omega_j t} \right)^{1/2}\int_{\R^n}e^{iS(t,x,y)/\eps} f(y)dy,
\end{equation*}
where
$$S(t,x,y) =  \sum_{j=1}^n
\frac{\omega_j}{\sin\omega_j t} \left( \frac{x_j^2+y_j^2}{2}\cos
\omega_j t -x_j y_j \right).$$
It is not hard to see that one can mimic the proof of
Th.~\ref{th:princ} to get the following,
\begin{theo}\label{th:anis}
Assume $2\leq n\leq 5$, $\frac{1}{2}<\sigma<\frac{2}{n-2}$, and
let $2 < r< \frac{2n}{n-2}$. If
$n=2$, there
exists $\delta>0$ such that in either of the two 
cases,
\begin{itemize}
\item $\sigma>\sigma_0(2)$, or
\item $\|f\|_\Sigma\leq \delta$,
\end{itemize}
the following holds (if $3\leq n \leq 5$, no additional assumption is
needed). Denote $\psi_\pm = W^{-1}_\pm f$, and
$$\varphi(t,x)=\frac{1}{2}\sum_{j=1}^n \frac{\omega_j
x_j^2}{\tan (\omega_j t)}.$$
Let $u^\eps$ be the solution of (\ref{eq:anis2}). 
Then  for $|t|<\frac{\pi}{\omega_d}$ (that is,
before refocusing), and
in $L^2\cap L^r$,
\begin{itemize}
\item If $0<t<\frac{\pi}{\omega_d}$, 
$$u^\eps(t,x)\Eq \eps 0 \prod_{j=1}^n\left(\frac{\omega_j}{2i\pi \sin
\omega_j t}\right)^{1/2}\widehat{\psi_+}\left(\frac{\omega_1 x_1}{\sin
\omega_1 t }, \ldots, \frac{\omega_n x_n}{\sin
\omega_n t } \right)e^{i\varphi(t,x)/\eps}.$$ 
\item If $-\frac{\pi}{\omega_d}<t<0$, 
$$u^\eps(t,x)\Eq \eps 0 \prod_{j=1}^n\left(\frac{\omega_j}{2i\pi \sin
\omega_j t}\right)^{1/2}\widehat{\psi_-}\left(\frac{\omega_1 x_1}{\sin
\omega_1 t }, \ldots, \frac{\omega_n x_n}{\sin
\omega_n t } \right)e^{i\varphi(t,x)/\eps}.$$
\end{itemize}
\end{theo}
At time $t=\frac{\pi}{\omega_d}$, the solution focuses on $E_d$,
and the nonlinear term is negligible (use the operators $H_j^\eps$ for
all indexes $j$ such that $\omega_j=\omega_d$, and $J^\eps_k$ for the
others). The nonlinear term will be
relevant again only if there exists a time where the focusings caused
by the different $\omega_j$'s ($1\leq j\leq d$) occur simultaneously,
that is if there are positive integers $\kappa_1,\ldots, \kappa_d$
such that 
$$t_1= \frac{\kappa_1\pi}{\omega_1}= \ldots =
\frac{\kappa_d\pi}{\omega_d}.$$
This means that the $\omega_j$'s are pairwise rationally
related. Therefore, at time $t= t_1$, the caustic crossing will be
described again by the scattering operator. Notice that when $t$
approaches $t_1$, the asymptotics given in Th.~\ref{th:anis} has been
modified in terms of Maslov indexes (for instance, since the crossing
of $E_d$ is linear, only linear phenomenon occur at leading order,
that is precisely a phase shift measured by the Maslov index). 
More precisely, for $(\kappa_d -1)\pi/\omega_d < t< t_1$, every
$\omega_j$ ($1\leq j\leq d$) part of the harmonic potential has caused
$\kappa_j -1$ (linear) caustic crossings, and 
$$u^\eps(t,x)\Eq \eps 0 \prod_{j=1}^d\left(\frac{\omega_j
e^{-i(\kappa_j-1)
\pi}}{2i\pi |\sin
\omega_j t|}\right)^{i_j/2}
\widehat{\psi_+}\left(\frac{\omega_1 x_1}{\sin
\omega_1 t }, \ldots, \frac{\omega_n x_n}{\sin
\omega_n t } \right)e^{i\varphi(t,x)/\eps}.$$
For $t_1<t<(\kappa_d +1)\pi/\omega_d $, one has,
$$u^\eps(t,x)\Eq \eps 0 \prod_{j=1}^d\left(\frac{\omega_j
e^{-i\kappa_j
\pi}}{2i\pi |\sin
\omega_j t|}\right)^{i_j/2}
\widehat{S\psi_+}\left(\frac{\omega_1 x_1}{\sin
\omega_1 t }, \ldots, \frac{\omega_n x_n}{\sin
\omega_n t } \right)e^{i\varphi(t,x)/\eps},$$
and so on. 

If the $\omega_j$'s are not pairwise rationally
related, then only linear phenomena occur near caustics, and they are
measured by Maslov indexes.

\providecommand{\bysame}{\leavevmode\hbox to3em{\hrulefill}\thinspace}
\providecommand{\href}[2]{#2}

\end{document}